\documentclass[aos,preprint]{imsart}

\usepackage{amsthm, amsmath, amssymb, amsfonts, graphicx, caption,subcaption,color}

\usepackage{tikz}
\usetikzlibrary{arrows}
\usepackage{caption}
\RequirePackage{algorithmic}
\RequirePackage[section]{algorithm}
\usepackage{xr-hyper}
\RequirePackage[square,authoryear]{natbib}
\usepackage{booktabs}
\usepackage{pifont}

\newcommand{\cmark}{\ding{51}}%
\newcommand{\xmark}{\ding{55}}%

\newtheorem{lemma}{Lemma}[section]
\newtheorem{theorem}{Theorem}[section]
\newtheorem{corollary}{Corollary}[section]
\newtheorem{remark}{Remark}[section]
\newtheorem{definition}{Definition}[section]
\newtheorem{notation}{Notation}[section]
\newtheorem{example}{Example}
\newtheorem*{assumption*}{Assumption}

\newcommand{\Prob}{\mathbb{P}}

\DeclareMathOperator{\Exp}{E}
\DeclareMathOperator{\Cov}{Cov}
\DeclareMathOperator{\Var}{Var}
\DeclareMathOperator{\CPDAG}{CPDAG}

\DeclareMathOperator{\skeleton}{skeleton}
\DeclareMathOperator{\NPN}{NPN}

\newcommand*{\indep}{\!\perp\!\!\!\perp}
\newcommand*{\ndsep}{\not\perp}
\newcommand*{\dsep}{\perp}

\begin{document}

\begin{frontmatter}

\title{High-dimensional consistency in score-based and hybrid structure learning}
\runtitle{High-dimensional consistency in structure learning}
\thankstext{T1}{Supported in part by Swiss NSF Grant 200021\_149760.}
\begin{aug}
\author{\fnms{Preetam} \snm{Nandy}\thanksref{m1, T1}\ead[label=e1]{preetamnandy@gmail.com},}
\author{\fnms{Alain} \snm{Hauser}\thanksref{m2}\ead[label=e2]{alhauser@google.com}}
\and
\author{\fnms{Marloes H.} \snm{Maathuis}\thanksref{m1, T1}\ead[label=e3]{maathuis@stat.math.ethz.ch}}

\affiliation{ETH Z\"urich\thanksmark{m1} and Bern University of Applied Sciences\thanksmark{m2}}

\runauthor{P. Nandy, A. Hauser and M. H. Maathuis}

\address{P. Nandy\\
        M. H. Maathuis\\
        ETH Z\"urich\\
        Seminar for Statistics\\
        R\"amistrasse 101\\
        8092 Z\"urich, Switzerland\\
        \printead{e1}\\
        \printead{e3}}

\address{A. Hauser\\
        Department of Engineering and \\
        Information Technology\\
        Bern University of Applied Sciences \\
        Jlcoweg 1 \\
        3400 Burgdorf, Switzerland\\
        (now at Google)\\
        \printead{e2}}
\end{aug}

\begin{abstract}
 Main approaches for learning Bayesian networks can be classified as constraint-based, score-based or hybrid methods. Although high-dimensional consistency results are available for constraint-based methods like the PC algorithm, such results have not been proved for score-based or hybrid methods, and most of the hybrid methods have not even shown to be consistent in the classical setting where the number of variables remains fixed and the sample size tends to infinity. In this paper, we show that consistency of hybrid methods based on greedy equivalence search (GES) can be achieved in the classical setting with adaptive restrictions on the search space that depend on the current state of the algorithm. Moreover, we prove consistency of GES and adaptively restricted GES (ARGES) in several sparse high-dimensional settings. ARGES scales well to sparse graphs with thousands of variables and our simulation study indicates that both GES and ARGES generally outperform the PC algorithm. 
\end{abstract}

\begin{keyword}
 \kwd{Bayesian network}
  \kwd{directed acyclic graph (DAG)}
 \kwd{linear structural equation model (linear SEM)}
 \kwd{structure learning}
 \kwd{greedy equivalence search (GES)}
 \kwd{score-based method}
 \kwd{hybrid method}
 \kwd{high-dimensional data}
 \kwd{consistency}
\end{keyword}

\end{frontmatter}
\section{Introduction}\label{section: introduction}
A Bayesian network consists of a directed acyclic graph (DAG) on a set of variables and conditional distributions for each node given its parents in the DAG. Bayesian networks can be used for various purposes, such as efficiently modeling the joint distribution of the variables, constructing decision support systems, probabilistic reasoning in expert systems, and causal inference. 

In a Bayesian network, the DAG encodes conditional independence relationships that must hold among the corresponding random variables. Several DAGs can encode exactly the same set of conditional independence relationships. Such DAGs are called Markov equivalent and form a Markov equivalence class (see Section \ref{subsection: joint distribution}). A Markov equivalence class can be uniquely represented by a completed partially directed acyclic graph (CPDAG). We consider estimating the CPDAG of a Bayesian network from observational data, and we refer to this as structure learning. Main approaches for structure learning can be classified as constraint-based, score-based or hybrid methods.

Constraint-based methods, such as the PC algorithm \citep{SpirtesEtAl00}, are based on conditional independence tests. The PC algorithm and its variants (e.g., \cite{HarrisDrton13, ColomboMaathuis14}) have been widely applied to high-dimensional datasets (e.g., \cite{MaathuisColomboKalischBuehlmann10, SchmidbergerEtAl11, StekhovenEtAl12, VerdugoEtAl13, LeEtAl13, GaoCui15}), partly because they were shown to be consistent in sparse high-dimensional settings where the number of variables is allowed to grow with the sample size \citep{KalischBuehlmann07a, HarrisDrton13, ColomboMaathuis14}, and partly because they scale well to sparse graphs with thousands of variables.

Score-based methods aim to optimize a scoring criterion over the space of possible CPDAGs or DAGs, typically through a greedy search procedure (see Section \ref{subsection: hill-climbing and GES}). Greedy equivalence search (GES) \citep{Chickering03} is a popular score-based method, which was shown to be consistent in the classical setting where the number of variables remains fixed and the sample size goes to infinity. This consistency result of \cite{Chickering03} is remarkable, since it involves a \emph{greedy} search. However, GES has not been shown to be consistent in high-dimensional settings. \cite{VandeGeerBuhlmann13} proved high-dimensional consistency of the global optimum of an $\ell_0$-penalized likelihood score function under the multivariate Gaussian assumption, but it has not been proved that a greedy search method like GES can find the global optimum in a high-dimensional setting. Another obstacle for applying score-based methods like GES to high-dimensional data is that they do not scale well to large graphs.

A hybrid method combines a score-based method either with a constraint-based method or with a variable selection method. Such methods often use a greedy search on a restricted search space in order to achieve computational efficiency, where the restricted space is estimated using conditional independence tests or variable selection methods \citep{TsamardinosEtAl06, SchmidtEtAl07, SchulteEtAL10, Alonso-BarbaEtAl13}. Common choices for the restricted search space are an estimated skeleton of the CPDAG (CPDAG-skeleton) or an estimated conditional independence graph (CIG). The CIG (or moral graph or Markov network) of a joint distribution of $p$ variables $X_1,\ldots,X_p$ is an undirected graph where two nodes $X_i$ and $X_j$ are adjacent if and only if they are conditionally dependent given $\{X_1,\ldots,X_p\}\setminus\{X_i,X_j\}$. The CIG is a supergraph of the CPDAG-skeleton. Hybrid algorithms generally scale well with respect to the number of variables, but their consistency results are generally lacking even in the classical setting, except for \cite{Alonso-BarbaEtAl13}.


In a preliminary simulation study, we compared the performances of PC, GES, and GES restricted to an estimated CIG (RGES-CIG) in high-dimensional settings (see Section \ref{A2-section: preliminary simulations} of the supplementary material). Table \ref{table: summary table 1} summarizes our findings from these preliminary simulations and the existence of consistency results in the literature.
\begin{table}[!ht]
  \caption{Summary of performance and existing consistency results, where tick marks represent good performance or existence of consistency results, cross marks represent bad performance, and question marks represent non-existence of consistency results.}
  \label{table: summary table 1}
  \begin{tabular}{lcccc}
    \toprule
                            & speed  & \begin{tabular}[b]{c} estimation\\performance \end{tabular} & consistency & \begin{tabular}[b]{c}high-dimensional\\ consistency \end{tabular}\\
    \midrule
    PC                    & \cmark & \xmark  &  \cmark & \cmark \\
    GES                   & \xmark & \cmark  & \cmark & ?  \\
    RGES-CIG             &  \cmark & \cmark  &  ?  & ? 
  \end{tabular}
\end{table}

Although GES and RGES-CIG outperform PC in terms of estimation performance in our high-dimensional simulations, we find that PC is the most commonly applied method in high-dimensional applications. We suspect that the main reasons for the lack of popularity of score-based and hybrid methods in high-dimensional applications are that they lack consistency results in high-dimensional settings and/or that they do not scale well to large graphs. In this paper, 
we prove high-dimensional consistency of GES, and we propose new hybrid algorithms based on GES that are consistent in several sparse high-dimensional settings \emph{and} scale well to large sparse graphs. 
To the best of our knowledge, these are the first results on high-dimensional consistency of score-based and hybrid methods.\\

The paper is organized as follows. Section \ref{section: preliminaries} discusses some necessary background  knowledge. In Section \ref{section: motivating example}, we show with an explicit example that naive hybrid versions of GES that restrict the search space to the CIG or the CPDAG-skeleton are inconsistent. 
This shows that the search path of GES may have to leave the search space determined by the CIG or the CPDAG-skeleton, even though the true CPDAG lies within these search spaces.

In Section \ref{section: ARGES} we provide a novel insight into how consistency can be achieved with hybrid algorithms based on GES, by imposing a restriction on the search space that changes adaptively depending on the current state of an algorithm. In particular, we propose a new method called adaptively restricted greedy equivalence search (ARGES), where in addition to the edges of the CIG (or the CPDAG-skeleton), we allow the shields of v-structures (or unshielded triples) in the current CPDAG, at every step of the algorithm. Our consistency proofs are based on a new characterization of independence maps (Theorem \ref{theorem: characterization of independence maps}), which is an interesting result in itself. 

In Section \ref{section: high-dimensional consistency} we prove consistency of GES and ARGES in certain sparse high-dimensional settings with multivariate Gaussian distributions. 
As a key ingredient of our proof, we show a connection between constraint-based and score-based methods. This connection enables us to extend our high-dimensional consistency results to linear structural equation models with sub-Gaussian errors (Section \ref{section: high dimensional consistency for LSEM}). Furthermore, it motivated us to define a scoring criterion based on rank correlations, and hence to derive high-dimensional consistency results for GES and ARGES for nonparanormal distributions (Section \ref{section: high dimensional consistency in the NPN setting}). This result is analogous to the high-dimensional consistency result of the Rank-PC algorithm \citep{HarrisDrton13}.  
Section \ref{section: high-dimensional simulations} contains simulation results, where we compare the finite sample performances and runtimes of PC, GES, ARGES and max-min hill-climbing \citep{TsamardinosEtAl06} in certain sparse high-dimensional settings. We end with a discussion and problems for future research in Section \ref{section: discussion}. 

All proofs are given in the supplementary material \citep{NandyHauserMaathuis15b}. An implementation of ARGES has been added to the \texttt{R}-package \textbf{pcalg} \citep{KalischEtAl12}.

\section{Preliminaries}\label{section: preliminaries}
\subsection{Graph terminology}\label{subsection: terminology}

We consider graphs $\mathcal{A} = (\mathbf{X},E)$, where the nodes (or vertices) $\mathbf{X} = \{X_1,\ldots,X_p\}$ represent random variables and the edges represent relationships between pairs of variables. The edges can be either \emph{directed} ($X_i \to X_k$) or \emph{undirected} ($X_i - X_k$). An \emph{(un)directed graph} can only contain (un)directed edges, whereas a \emph{partially directed graph} may contain both directed and undirected edges. For partially directed graphs $\mathcal{A} = (\mathbf{X},E)$ and $\mathcal{A}' = (\mathbf{X},E')$, we write $\mathcal{A} \subseteq \mathcal{A}'$ and $\mathcal{A} = \mathcal{A}'$ to denote $E \subseteq E'$ and $E = E'$ respectively. The \emph{skeleton} of a partially directed graph $\mathcal{A}$, denoted as $\skeleton(\mathcal{A})$, is the undirected graph that results from replacing all directed edges of $\mathcal{A}$ by undirected edges.

Two nodes $X_i$ and $X_k$ are \emph{adjacent} if there is an edge between them. Otherwise they are \emph{non-adjacent}. The set of all adjacent node of $X_i$ in $\mathcal{A}$ is denoted by $\mathbf{Adj}_{\mathcal{A}}(X_i)$. The degree of a node $X_i$ in $\mathcal{A}$ equals $|\mathbf{Adj}_{\mathcal{A}}(X_i)|$. If $X_i \to X_k$, then $X_i$ is a \emph{parent} of $X_k$. The set of all parents of $X_k$ in $\mathcal{A}$ is denoted by $\mathbf{Pa}_{\mathcal{A}}(X_k)$. A triple of nodes $(X_i,X_j,X_k)$ is an unshielded triple in $\mathcal{A}$ if $X_i$ and $X_k$ are non-adjacent in $\mathcal{A}$ and $\{X_i,X_k\} \subseteq \mathbf{Adj}_{\mathcal{A}}(X_j)$. An unshielded triple $(X_i,X_j,X_k)$ is a \emph{v-structure} if $X_i \to X_j \leftarrow X_k$. If $(X_i,X_j,X_k)$ is an unshielded triple in $\mathcal{A}$, then the edge $X_i \to X_k$ (or $X_k \to X_i$), which is not present in $\mathcal{A}$, is called a shield of the unshielded triple $(X_i,X_j,X_k)$.

A \emph{path} between $X_i$ and $X_k$ in a graph $\mathcal{A}$ is sequence of distinct nodes $(X_i,\dots,X_k)$ such that all pairs of successive nodes in the sequence are adjacent in $\mathcal{A}$. We use the shorthand $\pi_{\mathcal{A}}(X_i,\ldots,X_k)$ to denote a path in $\mathcal{A}$ with endpoint nodes $X_i$ and $X_k$. A non-endpoint node $X_{r}$ on a path $\pi_{\mathcal{A}}(X_i,\dots,X_{r-1},X_r,X_{r+1},\dots,X_k)$ is a \emph{collider} on the path if $X_{r-1} \to X_r \leftarrow X_{r+1}$. Otherwise it is a \emph{non-collider} on the path. An endpoint node on a path is neither a collider nor a non-collider on the path. A path without colliders is a \emph{non-collider path}. A path of two nodes is a trivial non-collider path. 

A \emph{directed path} from $X_i$ to $X_k$ is a path between $X_i$ and $X_k$, where all edges are directed towards $X_k$. If there is a directed path from $X_i$ to $X_k$, then $X_k$ is a \emph{descendant} of $X_i$. Otherwise it is a \emph{non-descendant}. We use the convention that each node is a descendant of itself. The set of all descendants (non-descendants) of $X_i$ in $\mathcal{A}$ is denoted by $\mathbf{De}_{\mathcal{A}}(X_i)$ ($\mathbf{Nd}_{\mathcal{A}}(X_i)$).

A path between $X_i$ and $X_k$ of at least three nodes, together with edge between $X_i$ and $X_k$ forms a \emph{cycle}. A directed path from $X_i$ to $X_k$ together with $X_k \to X_i$ forms a \emph{directed cycle}. A directed graph or partially directed graph without directed cycles is called \emph{acyclic}.  A graph that is both (partially) directed and acyclic, is a \emph{(partially) directed acyclic graph} or (P)DAG.

We will typically denote an arbitrary DAG by $\mathcal{H}$ (or $\mathcal{H}_n$), and an arbitrary partially directed graph by $\mathcal{A}$ (or $\mathcal{A}_n$). Graphs are always assumed to have vertex set $\mathbf{X} = \{X_1,\ldots,X_p\}$ (or $\mathbf{X}_n = \{X_{n1},\ldots,X_{np_n}\}$).


\subsection{Bayesian network terminology}\label{subsection: joint distribution}

We consider a random vector $\mathbf{X} = (X_1,\ldots,X_p)^T$ with a parametric density $f(\cdot)$.
The density \emph{factorizes} according to a DAG $\mathcal{H}$ if there exists a set of parameter values $\boldsymbol\Theta = \{\boldsymbol\theta_1,\ldots,\boldsymbol\theta_p\}$ such that
\vspace{-0.1in}
\begin{equation}\label{eq: factorization}
f(x_1,\ldots,x_p) = \prod_{i=1}^{p} f_i(x_i \mid \mathbf{Pa}_{\mathcal{H}}(X_i) = \mathbf{Pa}_{\mathcal{H}}(x_i), \boldsymbol\theta_i),
\vspace{-0.1in}
\end{equation}
where $\boldsymbol\theta_i$ specifies the conditional density of $X_i$ given its parents in $\mathcal{H}$ and $\mathbf{Pa}_{\mathcal{H}}(x_i)$ denotes the sub-vector of $(x_1,\ldots,x_p)$ that corresponds to $\mathbf{Pa}_{\mathcal{H}}(X_i)$. 
Such a pair $(\mathcal{H},\boldsymbol\Theta)$ is a Bayesian network that defines the joint distribution. The DAG $\mathcal{H}$ of a Bayesian network $(\mathcal{H},\boldsymbol\Theta)$ encodes conditional independence constraints that must hold in any distribution that factorizes according to $\mathcal{H}$. Conditional independence constraints encoded by a DAG can be read off from the DAG using the notion of \emph{d-separation}.

\begin{definition} (d-separation, see Definition 1.2.3 of \citet{PearlBook09})
   Let $\mathbf{S}$ be a subset of nodes in a DAG $\mathcal{H}$, where $\mathbf{S}$ does not contain $X_i$ and $X_k$. Then $\mathbf{S}$ \emph{blocks} a path $\pi_{\mathcal{H}}(X_i,\ldots,X_k)$ if at least one of the following holds: (i) $\pi_{\mathcal{H}}$ contains a non-collider that is in $\bf S$, or (ii) $\pi_{\mathcal{H}}$ contains a collider that has no descendant in $\bf S$. Otherwise $\pi_{\mathcal{H}}$ is \emph{open} given $\bf S$. 
   For pairwise disjoint sets of nodes $\mathbf{W}_1$, $\mathbf{W}_2$ and $\mathbf{S}$, we say that $\mathbf{W}_1$ and $\mathbf{W}_2$ are \emph{d-separated} by $\bf S$ in $\mathcal{H}$ if every path  between a node in $\mathbf{W}_1$ and a node in $\mathbf{W}_2$ is blocked by $\bf S$. This is denoted by $\mathbf{W}_1 \dsep_{\mathcal{H}} \mathbf{W}_2\mid \mathbf{S}$.  Otherwise, $\mathbf{W}_1$ and $\mathbf{W}_2$ are \emph{d-connected} given $\mathbf{S}$ in $\mathcal{H}$, denoted by $\mathbf{W}_1 \ndsep_{\mathcal{H}} \mathbf{W}_2\mid \mathbf{S}$.
\end{definition}

The distribution of $\mathbf{X}$ is \emph{DAG-perfect} \citep{Chickering03} if there exists a DAG $\mathcal{G}_0$ such that (i) every independence constraint encoded by $\mathcal{G}_0$ holds in the distribution of $\mathbf{X}$, and (ii) every independence constraint that holds in the distribution is encoded by $\mathcal{G}_0$. Such a DAG $\mathcal{G}_0$ is called a \emph{perfect map} of the distribution. Condition (i) is known as the \emph{global Markov property}, and condition (ii) is the so-called \emph{faithfulness} condition (see, for example, Definition 3.8 of \cite{KollerFriedman09}). In this paper, we only consider DAG-perfect distributions (as in \cite{Chickering03}).

DAGs that encode exactly the same set of conditional independence constraints form a \emph{Markov equivalence class} \citep{VermaPearl90}. Two DAGs belong to the same Markov equivalence class if and only if they have the same skeleton and the same v-structures \citep{VermaPearl90}. A Markov equivalence class of DAGs can be uniquely represented by a \emph{completed partially directed acyclic graph} (CPDAG), which is a PDAG that satisfies the following: $X_i\to X_k$ in the CPDAG if $X_i\to X_k$ in every DAG in the Markov equivalence class, and $X_i - X_k$ in the CPDAG if the Markov equivalence class contains a DAG in which $X_i\to X_k$ as well as a DAG in which $X_i\leftarrow X_k$ \citep{VermaPearl90, AnderssonEtAl97, Chickering03}. We use the notation $\CPDAG(\mathcal{H})$ to denote the CPDAG that represents the Markov equivalence class of a DAG $\mathcal{H}$.

\begin{notation}\label{notation: true graphs}  We reserve the notation $\mathcal{G}_0$, $\mathcal{C}_0$, $\mathcal{I}_0$ and $\mathcal{D}_n$ for the following:
  $\mathcal{G}_0$ denotes a perfect map of the distribution of $\mathbf{X}$ with $\mathcal{C}_0 = \CPDAG(\mathcal{G}_0)$, 
  $\mathcal{I}_0$ is the conditional independence graph of $\mathbf{X}$,
  and $\mathcal D_n$ denotes the data, consisting of $n$ i.i.d.\ observations of $\mathbf{X}$.
\end{notation}

  For a DAG $\mathcal{H}$, let $\mathbb{CI}(\mathcal{H})$ denote the set of all conditional independence constraints encoded by $\mathcal{H}$. By definition, $\mathcal{H}_1$ and $\mathcal{H}_2$ are Markov equivalent DAGs if and only if $\mathbb{CI}(\mathcal{H}_1) = \mathbb{CI}(\mathcal{H}_2)$. Thus, for a CPDAG $\mathcal{C}$, we use the notation $\mathbb{CI}(\mathcal{C})$ unambiguously to denote the set of all conditional independence constraints encoded by any DAG $\mathcal{H}$ in the Markov equivalence class of $\mathcal{C}$.

  \begin{definition}\label{definition: independence map}
  A (CP)DAG $\mathcal{A}$ is an \emph{independence map} of a (CP)DAG $\mathcal{A}'$ if $ \mathbb{CI}(\mathcal{A}) \subseteq \mathbb{CI}(\mathcal{A}')$.
  \end{definition}


\subsection{Properties of a scoring criterion}\label{subsection: properties scores}

Score-based and hybrid structure learning methods require a scoring criterion $\mathcal S(\mathcal H, \mathcal D_n)$ that measures the quality of a candidate DAG $\mathcal H$ with respect to given data $\mathcal D_n$. Throughout this paper, we assume without loss of generality that optimizing a scoring criterion corresponds to minimizing it. Hence, we say that the score improves by moving from $\mathcal{H}$ to $\mathcal{H}'$ if $\mathcal S(\mathcal H', \mathcal D_n) < \mathcal S(\mathcal H, \mathcal D_n)$.

We consider scoring criterions that are \emph{score equivalent}, \emph{decomposable} and \emph{consistent} (see Section \ref{A2-section: scoring criterion} of the supplementary materials). These properties are also assumed by \cite{Chickering03} as basic requirements of a scoring criterion to be used in GES. Score equivalence ensures that all DAGs in a Markov equivalence class get the same score, and the common score is defined as the score of the Markov equivalence class or its representative CPDAG. Decomposability of a scoring criterion facilitates fast computation of the score difference between two DAGs that differ by a few edges. Consistency of $\mathcal{S}$ assures that $\mathcal{G}_0$ has a lower score than any DAG that is not in the Markov equivalence class of $\mathcal{G}_0$, with probability approaching one as $n \to \infty$ (Proposition 8 of \cite{Chickering03}). For multivariate Gaussian distributions, the Bayesian information criterion (BIC) is an example of a scoring criterion that is score equivalent, decomposable and consistent. BIC was chosen as the scoring criterion of GES \citep{Chickering03}.


\subsection{Greedy equivalence search (GES) \citep{Chickering03}}\label{subsection: hill-climbing and GES}

GES is a greedy search algorithm that aims to optimize a score function on the space of CPDAGs. GES uses a forward phase and a backward phase. The forward phase starts with an initial CPDAG (often an empty graph) and sequentially obtains larger CPDAGs by adding exactly one edge at each step. Among all the possible single edge additions at every step, it selects the one that minimizes the score function. The forward phase ends when the score of the CPDAG can no longer be improved by a single edge addition. The backward phase starts with the output of the forward phase and sequentially obtains smaller CPDAGs by deleting exactly one edge at each step. It selects the optimal single edge deletion at each step and stops when the score can no longer be improved by a single edge deletion. 

Conceptually, the single edge additions and deletions are defined as follows. At every step in the forward (or backward) phase of GES, one can first list all DAGs in the Markov equivalence class of the current CPDAG, then consider all DAGs that can be obtained from the aforementioned DAGs by a single edge addition (or deletion), and finally move to the CPDAG that corresponds to a resulting DAG that has the minimum score (if the minimum score is smaller than the score of the current CPDAG). Thus, at each step, the skeleton of the CPDAG changes by exactly one edge, but the orientations may change for several edges. These moves, however, can be determined much more efficiently, and we refer to \cite{Chickering03} for details.

Pseudocodes of the forward and the backward phases are given in Section \ref{A2-section: pseudocode for GES} of the supplementary material and we refer to Figure \ref{A2-fig: stepwise (R)GES} of the supplementary material for an illustration of the search path of GES for Example \ref{example: inconsistency} (Section \ref{section: motivating example}).

\section{Inconsistency of GES restricted to the CIG/CPDAG-skeleton}\label{section: motivating example}

Naive hybrid versions of GES restrict the search space to an estimated CIG or CPDAG-skeleton. We refer to these hybrid methods as RGES-CIG or RGES-skeleton. More precisely, we restrict the search space of GES by allowing an edge $X_i \to X_j$ for addition only if $X_i$ and $X_j$ are adjacent in the CIG or in the CPDAG-skeleton. We will show inconsistency of these methods using an explicit example, where we assume that the CIG or the CPDAG-skeleton can be estimated consistently meaning that we restrict the search space to the true CIG or the true CPDAG-skeleton.

\begin{example}\label{example: inconsistency}
We consider $\mathbf{X} = (X_1,X_2,X_3,X_4)^T$ with a joint distribution defined via the following linear structural equation model (SEM): $X_1 = \epsilon_1$, $X_2 = \epsilon_2$, $X_3 = 1.4X_1 + 1.3X_2 + \epsilon_3$ and $X_4 = 1.2X_2 + 0.9X_3 + \epsilon_4$, where $\epsilon_1,\ldots,\epsilon_4$ are independently distributed standard Gaussian random variables. We write the linear SEM in matrix notation as $\mathbf{X} = B \mathbf{X} + \boldsymbol\epsilon$, where $B$ is a lower triangular matrix of coefficients and $\boldsymbol\epsilon = (\epsilon_1,\ldots, \epsilon_4)^T$. Thus $\mathbf{X} = (\mathrm{I} - B)^{-1}\boldsymbol\epsilon$ and $\mathbf{X}$ has a zero-mean multivariate Gaussian distribution with covariance matrix $\Sigma_0 = (\mathrm{I} - B)^{-1}(\mathrm{I} - B)^{-T}$.

This linear SEM can be represented by the DAG $\mathcal{G}_0$ in Figure \ref{fig: example}(a), where an edge $X_i \to X_j$ is present if and only if $B_{ji} \neq 0$ and then the weight of the edge $X_i \to X_j$ is $B_{ji}$. Note that $\mathcal{G}_0$ is a perfect map of the distribution of $\mathbf{X}$. The Markov equivalence class of $\mathcal{G}_0$ contains only one DAG and hence the corresponding CPDAG is identical to $\mathcal{G}_0$. The CIG contains all possible undirected edges except for the edge $X_1 - X_4$ (Figure \ref{fig: example}(b)), whereas the CPDAG-skeleton additionally excludes the undirected edge $X_1 - X_2$ (Figure \ref{fig: example}(c)).

\begin{figure}[!t]
  \centering
  \footnotesize
  \begin{subfigure}[t]{0.25\textwidth}
    \centering
       \begin{tikzpicture}[scale=0.8, transform shape]
        \node         (X3) at (0,0)                      {$X_3$};
   \node         (X4) at (1.5,0)	{$X_4$};
   \node         (X1) at (0,1.25) 	{$X_1$};
   \node     	(X2) at (1.5,1.25)        	{$X_2$};
   \draw[->, thick] (X1) edge  node[left=2pt]{$1.4$}  (X3);
            \draw[->, thick] (X2) edge  node[left =4pt, above=2pt]{$1.3$} (X3);
      \draw[->, thick] (X2) edge   node[right=2pt]{$1.2$} (X4);
   	\draw[->, thick] (X3) edge  node[below=3pt]{$0.9$}  (X4);
   \end{tikzpicture}
   \caption{}
   \label{subfig: WDAG}
\end{subfigure}
       \begin{subfigure}[t]{0.25\textwidth}
           \centering
       \begin{tikzpicture}[scale=0.8, transform shape]
       \centering
         \node         (X3) at (0,0)                      {$X_3$};
   \node         (X4) at (1.5,0)	{$X_4$};
   \node         (X1) at (0,1.25) 	{$X_1$};
   \node     	(X2) at (1.5,1.25)        	{$X_2$};
   \draw[-, thick] (X1) edge   (X3);
      \draw[-, thick] (X1) edge   (X2);
   \draw[-, thick] (X4) edge node[below=8.5pt]{}   (X3);
   \draw[-, thick] (X4) edge   (X2);
   \draw[-, thick] (X2) edge  (X3);
   \end{tikzpicture}
     \caption{}
     \label{subfig: CIG}
     \end{subfigure}
  \begin{subfigure}[t]{0.25\textwidth}
    \centering
       \begin{tikzpicture}[scale=0.8, transform shape]
        \node         (X3) at (0,0)                      {$X_3$};
   \node         (X4) at (1.5,0)	{$X_4$};
   \node         (X1) at (0,1.25) 	{$X_1$};
   \node     	(X2) at (1.5,1.25)        	{$X_2$};
   \draw[-, thick] (X1) edge   (X3);
            \draw[-, thick] (X2) edge   (X3);
      \draw[-, thick] (X2) edge   (X4);
   	\draw[-, thick] (X3) edge node[below=8.5pt]{}  (X4);
   \end{tikzpicture}
   \caption{}
   \label{subfig: CPDAG}
\end{subfigure}
  \caption{The weighted DAG in (a) represents the data generating process. The corresponding CIG and the CPDAG-skeleton are given in (b) and (c) respectively.}
 \label{fig: example}
  \end{figure}
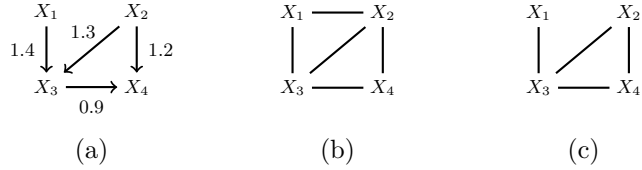

\begin{figure}[!ht]
  \centering
  \begin{subfigure}[t]{0.25\textwidth}
    \centering
       \begin{tikzpicture}[scale=0.8, transform shape]
        \node         (X3) at (0,0)                      {$X_3$};
   \node         (X4) at (1.5,0)	{$X_4$};
   \node         (X1) at (0,1.25) 	{$X_1$};
   \node     	(X2) at (1.5,1.25)        	{$X_2$};
   \draw[->, thick] (X1) edge  (X3);
            \draw[->, thick] (X2) edge   (X3);
      \draw[->, thick] (X2) edge  (X4);
   	\draw[->, thick] (X3) edge  (X4);
   \end{tikzpicture}
   \caption{}
   \label{subfig: GES}
\end{subfigure}
  \centering
  \begin{subfigure}[t]{0.25\textwidth}
    \centering
       \begin{tikzpicture}[scale=0.8, transform shape]
         \node         (X3) at (0,0)                      {$X_3$};
   \node         (X4) at (1.5,0)	{$X_4$};
   \node         (X1) at (0,1.25) 	{$X_1$};
   \node     	(X2) at (1.5,1.25)        	{$X_2$};
   \draw[->, thick] (X1) edge   (X3);
      \draw[->, thick] (X1) edge   (X2);
   \draw[->, thick] (X4) edge   (X3);
   \draw[->, thick] (X4) edge   (X2);
   \draw[-, thick] (X2) edge  (X3);
   \end{tikzpicture}
     \caption{}
     \label{subfig: RGES}
\end{subfigure}
  \begin{subfigure}[t]{0.25\textwidth}
    \centering
       \begin{tikzpicture}[scale=0.8, transform shape]
          \node         (X3) at (0,0)                      {$X_3$};
   \node         (X4) at (1.5,0)	{$X_4$};
   \node         (X1) at (0,1.25) 	{$X_1$};
   \node     	(X2) at (1.5,1.25)        	{$X_2$};
   \draw[->, thick] (X1) edge   (X3);
   \draw[->, thick] (X4) edge  (X3);
   \draw[->, thick] (X4) edge   (X2);
   \draw[->, thick] (X3) edge  (X2);
   \end{tikzpicture}
     \caption{}
     \label{subfig: RGES-skeleton}
     \end{subfigure}
\caption{The CPDAGs in (a), (b) and (c) are the large sample limit outputs of GES, RGES-CIG and RGES-skeleton respectively.} 
 \label{fig: large sample limits}
  \end{figure}
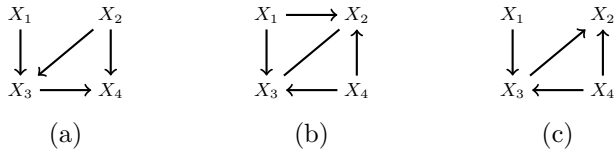

Figure \ref{fig: large sample limits} shows that the large sample limit output (with the BIC criterion) of GES equals $\mathcal{G}_0 = \CPDAG(\mathcal{G}_0)$, but the large sample limit outputs (with the BIC criterion) of RGES-CIG and RGES-skeleton are different from $\mathcal{G}_0$. The corresponding search paths of all three algorithms are given in Section \ref{A2-subsection: stepwise illustration} of the supplementary material. 

We initialized all algorithms by the empty graph. We determined the large sample limit outputs of all algorithms by choosing the scoring criterion to be the expected negative log-likelihood scores. A detailed description is given in Section \ref{A2-subsection: limit computation} of the supplementary material.

\end{example}


\begin{remark}\label{remark: inconsistency of some other hybrid algorithms}
We note that the linear Gaussian SEM given in Example \ref{example: inconsistency} can also be used to show inconsistency of hill-climbing DAG search, hill-climbing DAG search restricted to the CIG (e.g.\ \cite{SchmidtEtAl07}) and hill-climbing DAG search restricted to the CPDAG-skeleton (e.g.\ the max-min hill-climbing algorithm of \cite{TsamardinosEtAl06}). Details are given in Section \ref{A2-subsection: limit outputs} of the supplementary material.
\end{remark}

\begin{remark}
The consistency of score-based and hybrid algorithms corresponds to the soundness of these algorithms with an appropriate oracle scoring criterion. An oracle scoring criterion depends on the joint distribution instead of the data and can be viewed as the large sample limit of its finite sample counterpart. For example, the (penalized) expected log-likelihood score is an oracle score which is the large sample limit of the (penalized) average log-likelihood. In fact, Example \ref{example: inconsistency} shows that RGES-CIG and RGES-skeleton with the expected negative log-likelihood scoring criterion are not sound (cf.\ Section \ref{A2-subsection: limit computation} of the supplementary material).
\end{remark}

\section{Adaptively restricted greedy equivalence search (ARGES)}\label{section: ARGES}

In the previous section, we have seen that naive hybrid versions of GES can be inconsistent, although GES is consistent \citep{Chickering03}. We now propose adaptively restricted hybrid versions of GES that will be shown to be consistent in Section \ref{subsection: consistency of ARGES}.

We recall that the consistency proof of \citet{Chickering03} consists of two main steps. The first step shows that the output of the forward phase is an independence map of $\mathcal{C}_0$, with probability approaching one. 
The second step shows that $\Prob (\tilde{\mathcal{C}}_{n} = \mathcal{C}_0) \to 1$ given the result of the first step (Lemma 10 of \cite{Chickering03}), where $\tilde{\mathcal{C}}_{n}$ denotes the output of GES based on $n$ i.i.d.\ samples.

We consider hybrid versions of GES that modify only the forward phase of GES, by restricting edge additions. To retain consistency in such hybrid versions, it therefore suffices to ensure that the output of the forward phase is an independence map of $\mathcal{C}_0$ with probability approaching one. This motivated us to provide the following novel characterization of independence maps in Theorem \ref{theorem: characterization of independence maps}. Our adaptively restricted hybrid versions and their consistency will follow as a natural consequence of Theorem \ref{theorem: characterization of independence maps} and Chickering's consistency proof of GES.

\begin{theorem}\label{theorem: characterization of independence maps}
A DAG $\mathcal{H}$ is not an independence map of a DAG $\mathcal{G}$ if and only if 
\begin{enumerate}
\item $\skeleton(\mathcal{G}) \nsubseteq \skeleton(\mathcal{H})$, or 
\item there exists a triple of nodes $\{X_i,X_j,X_k\}$ such that $X_i$ and $X_k$ are non-adjacent in $\mathcal{H}$, $\pi_{\mathcal{H}}(X_i,X_j,X_k)$ is a non-collider path, and $\pi_{\mathcal{G}}(X_i,X_j,X_k)$ is a v-structure, or
\item there exists a triple of nodes $\{X_i,X_j,X_k\}$ such that $\pi_{\mathcal{H}}(X_i,X_j,X_k)$ is a v-structure and $X_i \ndsep_{\mathcal{G}} X_k \mid \mathbf{Pa}_{\mathcal{H}}(X_k)$, where without loss of generality we assume $X_i \in \mathbf{Nd}_{\mathcal{H}}(X_k)$.
\end{enumerate}
\end{theorem}

We note that Proposition 27 and Lemma 28 of \cite{Chickering03} imply that if one of the first two conditions of Theorem \ref{theorem: characterization of independence maps} hold, then $\mathcal{H}$ is not an independence map. If the third condition of Theorem \ref{theorem: characterization of independence maps} holds, then $X_i \ndsep_{\mathcal{G}} X_k \mid \mathbf{Pa}_{\mathcal{H}}(X_k)$ and $X_i \dsep_{\mathcal{H}} X_k \mid \mathbf{Pa}_{\mathcal{H}}(X_k)$ (since $X_i \in \mathbf{Nd}_{\mathcal{H}}(X_k) \setminus \mathbf{Pa}_{\mathcal{H}}(X_k)$), and hence $\mathcal{H}$ is not an independence map of $\mathcal{G}$. 

The proof of the ``only if'' part of theorem is rather involved and we provide some intuition by considering two special cases in Section \ref{A2-subsection: proof of the characterization theorem} of the supplementary material.\\

We will use Theorem \ref{theorem: characterization of independence maps} with $\mathcal{G} = \mathcal{G}_0$ to derive that consistency can be achieved with a hybrid algorithm that allows an edge addition between $X_i$ and $X_k$ at any step in the forward phase if (i) $X_i$ and $X_k$ are adjacent in $\skeleton(\mathcal{G}_0)$, (ii) $\pi_{\mathcal{G}_0}(X_i,X_j,X_k)$ is a v-structure for some $X_j$, or (iii) $(X_i,X_j,X_k)$ is a v-structure in the current CPDAG for some $X_j$.

Recall that $X_i$ and $X_k$ are adjacent in the CIG of the distribution of $\mathbf{X}$ if and only if $X_i$ and $X_k$ are adjacent in $\mathcal{G}_0$ or $\pi_{\mathcal{G}_0}(X_i,X_j,X_k)$ is a v-structure for some $X_j$. Thus, we modify RGES-CIG by additionally allowing edges that are shields of v-structures in the current CPDAG at every step of the forward phase.

RGES-skeleton allows an edge addition between $X_i$ and $X_k$ only if $X_i$ and $X_k$ are adjacent in $\skeleton(\mathcal{G}_0)$. Therefore, we modify RGES-skeleton by additionally allowing edges that are shields of unshielded triples in the current CPDAG at every step of the forward phase. 

We call these modified versions ARGES-CIG or ARGES-skeleton, and we describe them in detail below. We often refer to both ARGES-CIG and ARGES-skeleton as ARGES in statements that hold for both of them.

\subsection{The ARGES algorithm}\label{subsection: forward phase of ARGES}
Given an estimated CIG (or CPDAG-skeleton), ARGES greedily optimizes a scoring criterion in two phases: a forward phase that depends on the estimated CIG (or CPDAG-skeleton) and on the current state of the algorithm, and a backward phase that is identical to the backward phase of GES.

The pseudocode of the forward phase of ARGES is given in Algorithm \ref{algorithm: the forward phase of ARGES}. It starts with an initial CPDAG $\mathcal{C}_\text{start}$ (often an empty graph) and sequentially obtains a larger CPDAG by adding exactly one edge at each step. At every step, it selects an optimal move (that minimizes the score) from a given set of \emph{admissible moves}, which depend on an estimated CIG (or CPDAG-skeleton) and the current state of the algorithm. To define the set of admissible moves, we introduce the notion of \emph{admissible edge} for ARGES-CIG and ARGES-skeleton.

\begin{definition}\label{definition: admissible edge} (Admissible edge for ARGES-CIG)
   Let $\{X_i,X_k\}$ be a pair of non-adjacent nodes in a CPDAG $\mathcal{C}$. Then an edge between $X_i$ and $X_k$
   is \emph{admissible} for $\mathcal{C}$ with respect to an undirected graph $\mathcal{I}$ if at least one of the following hold:
   \begin{enumerate}
      \item $X_i$ and $X_k$ are adjacent in $\mathcal{I}$; or
      \item There exists a node $X_j$ such that $(X_i, X_j, X_k)$ is a v-structure in $\mathcal{C}$.
   \end{enumerate}
\end{definition}

If $\mathcal{I}$ equals the true CIG $\mathcal{I}_0$ or an estimate thereof, then the first condition of Definition \ref{definition: admissible edge} corresponds to the restriction of the search space to the (estimated) CIG. The second condition is our adaptive relaxation that allows shields of v-structures. 

\begin{definition}\label{definition: admissible edge for ARGES-skeleton} (Admissible edge for ARGES-skeleton)
   Let $\{X_i,X_k\}$ be a pair of non-adjacent nodes in a CPDAG $\mathcal{C}$. Then an edge between $X_i$ and $X_k$
   is \emph{admissible} for $\mathcal{C}$ with respect to an undirected graph $\mathcal{U}$ if at least one of the following hold:
   \begin{enumerate}
      \item $X_i$ and $X_k$ are adjacent in $\mathcal{U}$; or
      \item There exists a node $X_j$ such that $(X_i, X_j, X_k)$ is an unshielded triple in $\mathcal{C}$.
   \end{enumerate}
\end{definition}

For ARGES-skeleton, we will choose $\mathcal{U}$ to be (an estimate of) $\skeleton(\mathcal{C}_0)$, which is typically a smaller graph than (an estimate of) $\mathcal{I}_0$. Therefore, the first condition in Definition \ref{definition: admissible edge for ARGES-skeleton} is more restrictive than the first condition of Definition \ref{definition: admissible edge} when $\mathcal{U} = \skeleton(\mathcal{C}_0)$ and $\mathcal{I} = \mathcal{I}_0$. This is somewhat compensated by the second condition, which is less restrictive in Definition \ref{definition: admissible edge for ARGES-skeleton} than in Definition \ref{definition: admissible edge}.

\begin{definition}\label{definition: admissible move}(Admissible move)
   Let $\mathcal{A}$ be an undirected graph and $\mathcal{C}$ a CPDAG, such that the edge $X_i \to X_k$ is admissible for a CPDAG $\mathcal{C}$ with respect to $\mathcal{A}$, where we apply Definition \ref{definition: admissible edge} if $\mathcal{A}$ is (an estimate of) $\mathcal{I}_0$ and Definition \ref{definition: admissible edge for ARGES-skeleton} if $\mathcal{A}$ is (an estimate of) $\skeleton(\mathcal{C}_0)$. Then the move from $\mathcal{C}$ to another CPDAG $\mathcal{C}'$ is admissible with respect to $\mathcal{A}$ if there exist DAGs $\mathcal{H}$ and $\mathcal{H}'$ in the Markov equivalence classes described by $\mathcal{C}$ and $\mathcal{C}'$ respectively, such that $\mathcal{H}'$ can be obtained from $\mathcal{H}$ by adding the edge $X_i \to X_k$ (i.e., $X_i \in \mathbf{Nd}_{\mathcal{H}}(X_k) \setminus \mathbf{Pa}_{\mathcal{H}}(X_k)$).
\end{definition}


\begin{algorithm}[H]
   \caption{The forward phase of ARGES (based on $\mathcal{A}$)}\label{algorithm: the forward phase of ARGES}
   \begin{algorithmic}[1]\label{algorithm: ARGES}
   \REQUIRE A scoring criterion $\mathcal{S}$, the data $\mathcal{D}_n$, an initial CPDAG $\mathcal{C}_\text{start}$, and an undirected graph $\mathcal{A}$ (that equals either an (estimated) CIG or an (estimated) CPDAG-skeleton).
   \ENSURE A CPDAG
   \STATE $\mathcal{C}_\text{new} \leftarrow \mathcal{C}_\text{start}$;
   \REPEAT
   \STATE   $\hat{\mathcal{C}}^{f}_n \leftarrow \mathcal{C}_\text{new}$;
   \STATE $\mathfrak{C} \leftarrow$ the set of all CPDAGs $\mathcal{C}$ such that $\mathcal{S}(\mathcal{C},\mathcal{D}_n) < \mathcal{S}(\hat{\mathcal{C}}^{f}_n,\mathcal{D}_n)$, and $\mathcal{C}$ can be obtained by an admissible move from $\hat{\mathcal{C}}^{f}_n$ with respect to $\mathcal{A}$; \label{line: admissibility}
   \IF{$\mathfrak{C} \neq \emptyset$}
   \STATE choose $\mathcal{C}_\text{new}$ to be the CPDAG that minimizes the scoring criterion among the CPDAGs in $\mathfrak{C}$; \label{line: optimal choice 1}
   \ENDIF
   \UNTIL{$\mathfrak{C} = \emptyset$};
   \RETURN $\hat{\mathcal{C}}^{f}_n$.
   \end{algorithmic}
\end{algorithm}

The forward phase of ARGES resembles the forward phase of GES (Algorithm \ref{A2-algorithm: the forward phase of GES} of the supplementary material), with the difference that at each step an edge between two non-adjacent nodes can only be added if (i)
the nodes are adjacent in an estimated CIG (or CPDAG-skeleton), or (ii) the edge shields a v-structure (or an unshielded triple) in the current CPDAG. Therefore, the forward phase of GES is the same as Algorithm \ref{algorithm: ARGES} with $\mathcal{A}$ being the complete undirected graph. 

\subsection{Consistency in the classical setting}\label{subsection: consistency of ARGES}

In this subsection we prove consistency of ARGES in the classical setting, where the sample size $n$ tends to infinity and the number of variables $p$ remains fixed. We fix an initial CPDAG $\mathcal{C}_{start}$ and a score equivalent and consistent scoring criterion $\mathcal{S}$ (see Section \ref{subsection: properties scores}). 

As we discussed before, it suffices to show that the output of the forward phase of ARGES is an independence map of $\mathcal{C}_0$. In the proof of Lemma 9 of \cite{Chickering03}, Chickering argued that if the output of the forward phase of GES is not an independence map, its score can be improved (asymptotically) by an edge addition, which is a contradiction. The following corollary of Theorem \ref{theorem: characterization of independence maps} shows that if $\mathcal{H}$ is not an independence map of $\mathcal{G}_0$, then the score can be improved (asymptotically) by adding an admissible edge  $X_i \to X_k$ (as in Definition \ref{definition: admissible edge} or Definition \ref{definition: admissible edge for ARGES-skeleton}). This additional result allows us to follow Chickering's argument for showing that the output of the forward phase of ARGES is an independence map of $\mathcal{C}_0$, leading to the consistency results given in Theorems \ref{theorem: consistency of ARGES} and \ref{theorem: consistency of ARGES-skeleton}.

\begin{corollary}\label{corollary: improvable by an admissible edge}
If $\mathcal{H}$ is not an independence map of $\mathcal{G}_0$, then there exists a pair of non-adjacent nodes $\{X_i, X_k\}$ in $\mathcal{H}$ such that $X_i \in \mathbf{Nd}_{\mathcal{H}}(X_k)$, $X_i \ndsep_{\mathcal{G}_0} X_k \mid \mathbf{Pa}_{\mathcal{H}}(X_k)$, and the edge $X_i \to X_k$ is admissible for $\CPDAG(\mathcal{H})$ with respect to $\mathcal{I}_0$ (for ARGES-CIG) and with respect to $\skeleton(\mathcal{C}_0)$ (for ARGES-skeleton).
\end{corollary}


\begin{theorem}\label{theorem: consistency of ARGES}
Let the CPDAG $\hat{\mathcal C}_n$ be the output of ARGES-CIG based on the estimated CIG $\hat{\mathcal{I}}_n$, where $\hat{\mathcal{I}}_n$ satisfies
 $\underset{n \to \infty}{\lim}~\Prob (\hat{\mathcal{I}}_n \supseteq \mathcal{I}_0) =1.$ Then $\underset{n \to \infty}{\lim}~\Prob(\hat{\mathcal{C}}_n = \mathcal{C}_0) = 1$.
\end{theorem}


\begin{theorem}\label{theorem: consistency of ARGES-skeleton}
Let the CPDAG $\breve{\mathcal C}_n$ be the output of ARGES-skeleton based on the estimated CPDAG-skeleton $\hat{\mathcal{U}}_n$, where $\hat{\mathcal{U}}_n$ satisfies
 $\underset{n \to \infty}{\lim}~\Prob (\hat{\mathcal{U}}_n \supseteq \skeleton(\mathcal{C}_0)) =1.$ Then $\underset{n \to \infty}{\lim}~\Prob(\breve{\mathcal C}_n = \mathcal{C}_0) = 1$.
\end{theorem}

\section{High-dimensional consistency of GES and ARGES in the multivariate Gaussian setting}\label{section: high-dimensional consistency}


We prove high-dimensional consistency of (AR)GES with an $\ell_0$-penalized log-likelihood scoring criterion $\mathcal{S}_{\lambda}$ (Definition \ref{definition: penalized likelihood score}), using the following steps. We define a collection of oracle versions of (AR)GES with an oracle scoring criterion $\mathcal{S}_{\lambda}^*$ (Definition \ref{definition: oracle score}), and prove soundness of them. We complete the proof by showing that the sample version of (AR)GES with scoring criterion $\mathcal{S}_{\lambda_n}$ is identical to one of the oracle versions with probability approaching one, for a suitably chosen sequence of penalty parameters $\lambda_n$. 



\subsection{$\ell_0$-penalized log-likelihood score in the multivariate Gaussian setting}\label{subsection: l0 penalized likelihood score}

\begin{definition}\label{definition: penalized likelihood score}
Let $\mathcal{H} = (\mathbf{X},E)$ be a DAG. The $\ell_0$-penalized log-likelihood score with penalty parameter $\lambda$ is given by 
\vspace{-0.05in}
\begin{align*}
\mathcal{S}_{\lambda}(\mathcal{H},\mathcal{D}_n) = -  \sum_{i=1}^{p}\frac{1}{n} \log \left(L(\hat{\boldsymbol\theta}_{i}(\mathcal{H}),\mathcal{D}_n(X_i)|\mathcal{D}_n(\mathbf{Pa}_{\mathcal{H}}(X_i)))\right) + \lambda |E|,
\end{align*}
where $L({\boldsymbol\theta}_{i}(\mathcal{H}),\mathcal{D}_n(X_i)|\mathcal{D}_n(\mathbf{Pa}_{\mathcal{H}}(X_i)))$ is the likelihood function that corresponds to the conditional density of $X_i$ given $\mathbf{Pa}_{\mathcal{H}}(X_i)$ and $$\hat{\boldsymbol\theta}_{i}(\mathcal{H}) = \mathop{\mathrm{argmax}}_{{\boldsymbol\theta}_{i}(\mathcal{H})}~ \frac{1}{n} \log \Big(L({\boldsymbol\theta}_{i}(\mathcal{H}),\mathcal{D}_n(X_i)|\mathcal{D}_n(\mathbf{Pa}_{\mathcal{H}}(X_i))) \Big)$$ is the maximum likelihood estimate (MLE) of the parameter vector $\boldsymbol\theta_{i}(\mathcal{H})$.
\end{definition}

\begin{remark}\label{remark: BIC}
The BIC criterion is a special case of the $\ell_0$-penalized log-likelihood score. In particular, the BIC score of a DAG $\mathcal{H}$ equals $2n\mathcal{S}_{\lambda_n}(\mathcal{H},\mathcal{D}_n)$ with $\lambda_n = \frac{\log(n)}{2n}$.
\end{remark}

The following lemma shows that when the distribution of $\mathbf{X}$ is multivariate Gaussian and $\lambda$ is suitably chosen, the $\ell_0$-penalized log-likelihood score of a DAG $\mathcal{H}$ can be improved by adding an edge $X_i \rightarrow X_j$ if and only if the sample partial correlation between $X_i$ and $X_j$ given $ \mathbf{Pa}_{\mathcal{H}}(X_j)$ is nonzero. This is one of the key results for our proof of high-dimensional consistency of (AR)GES.

\begin{lemma}\label{lemma: score difference}
Let $\mathcal{H} = (\mathbf{X},E)$ be a DAG such that $X_i\in \mathbf{Nd}_{\mathcal{H}}(X_k)\setminus \mathbf{Pa}_{\mathcal{H}}(X_k)$. Let $\mathcal{H}'=(\mathbf{X},E\cup \{X_i \rightarrow X_k\})$. If the distribution of $\mathbf{X}$ is multivariate Gaussian, then the $\ell_0$-penalized log-likelihood score difference between $\mathcal{H}'$ and $\mathcal{H}$ is
\vspace{-0.05in}
\begin{align}\label{eq: score difference formula}
\mathcal{S}_{\lambda}(\mathcal{H}',\mathcal{D}_n) - \mathcal{S}_{\lambda}(\mathcal{H},\mathcal{D}_n) = \frac{1}{2}\log \left(1 - \hat{\rho}_{ik|\mathbf{Pa}_{\mathcal{H}}(k)}^2\right) + \lambda,
\end{align}
where $\hat{\rho}_{ik|\mathbf{Pa}_{\mathcal{H}}(k)}$ denotes the sample partial correlation between $X_i$ and $X_k$ given $ \mathbf{Pa}_{\mathcal{H}}(X_k)$.
\end{lemma}

The first term on the right-hand side of \eqref{eq: score difference formula} equals the negative of the conditional mutual information between $X_i$ and $X_k$ given $\mathbf{Pa}_{\mathcal{H}}(X_k)$. Thus, Lemma \ref{lemma: score difference} shows that score-based methods like GES essentially use conditional independence tests (cf.\ \cite{AnandkumarEtAl12}) for sequentially adding and deleting edges starting from an initial graph. This shows that the score-based GES algorithm and the constraint-based PC algorithm use the same basic principle (a conditional independence test) in the multivariate Gaussian setting. Although this connection (see Section \ref{section: discussion} for more details) between PC and GES is not very surprising, we were unable to find it in the literature.

Further, Lemma \ref{lemma: score difference} also opens the possibility to define generalized scoring criterions, by replacing the Gaussian conditional mutual information in (2) with a more general measure of conditional independence. In fact, we exploit this in Section 6 to extend our high-dimensional consistency results to nonparanormal distributions.


We define an oracle version of the $\ell_0$- penalized log-likelihood scoring criterion by replacing the log-likelihood in Definition \ref{definition: penalized likelihood score} by its expectation with respect to the distribution $F$ of $\mathbf{X}$. The oracle score given by Definition \ref{definition: oracle score} will be used to define a collection of oracle versions of (AR)GES in Section \ref{subsection: high-dimensional consistency of ARGES}.

\begin{definition}\label{definition: oracle score} (Oracle score)
Let $\mathcal{H} = (\mathbf{X},E)$ be a DAG. We define the oracle score of $\mathcal{H}$ with respect to the distribution $F$ of $\mathbf{X}$ as 
\vspace{-0.05in}
\begin{align*}
\mathcal{S}_{\lambda}^*(\mathcal{H},F) = -  \sum_{i=1}^{p} \Exp \left[ \log \left(L({\boldsymbol\theta}_{i}^*(\mathcal{H}),X_i|\mathbf{Pa}_{\mathcal{H}}(X_i))\right) \right] + \lambda |E|,
\end{align*}
where $L({\boldsymbol\theta}_{i},X_i|\mathbf{Pa}_{\mathcal{H}}(X_i))$ is the likelihood function that corresponds to the conditional density of $X_i$ given $\mathbf{Pa}_{\mathcal{H}}(X_i)$ and \vspace{-0.05in}$${\boldsymbol\theta}_{i}^*(\mathcal{H}) = \underset{{\boldsymbol\theta}_{i}}{\mathrm{argmax}}~  \Exp \left[ \log \left( L({\boldsymbol\theta}_{i},X_i|\mathbf{Pa}_{\mathcal{H}}(X_i))\right) \right].$$
\end{definition}


We note that both $\mathcal{S}_{\lambda}$ and $\mathcal{S}_{\lambda}^*$ are decomposable (see Definition \ref{A2-definition: decomposable} of the supplementary material). Moreover, they are score equivalent (see Definition \ref{A2-definition: score equivalent} of the supplementary material) when the distribution $F$ is multivariate Gaussian. The scoring criterion used to compute the large sample outputs in Example \ref{example: inconsistency} equals $\mathcal{S}_{\lambda}^*$ with $\lambda = 0$ (see Section \ref{A2-subsection: limit computation} of the supplementary material).

The following lemma is analogous to Lemma \ref{lemma: score difference}. We do not provide a proof of Lemma \ref{lemma: oracle score difference}, since it can be obtained from the proof of Lemma \ref{lemma: score difference} by replacing the sample quantities by the corresponding population quantities (e.g.\ averages should be replaced by expectations and sample regression coefficients should be replaced by their population counterparts).

\begin{lemma}\label{lemma: oracle score difference}
Let $\mathcal{H} = (\mathbf{X},E)$ be a DAG such that $X_i \in \mathbf{Nd}_{\mathcal{H}}(X_k)\setminus \mathbf{Pa}_{\mathcal{H}}(X_k)$. Let $\mathcal{H}'=(\mathbf{X},E\cup \{X_i \rightarrow X_k\})$. If the distribution $\mathbf{X}$ is multivariate Gaussian, then the oracle score difference between $\mathcal{H}'$ and $\mathcal{H}$ is
\begin{align*}
\mathcal{S}_{\lambda}^*(\mathcal{H}',F) - \mathcal{S}_{\lambda}^*(\mathcal{H},F) = \frac{1}{2}\log \left(1 - \rho_{ik|\mathbf{Pa}_{\mathcal{H}}(k)}^2\right) + \lambda,
\end{align*}
where $\rho_{ik|\mathbf{Pa}_{\mathcal{H}}(k)}$ denotes the partial correlation between $X_i$ and $X_k$ given $ \mathbf{Pa}_{\mathcal{H}}(X_k)$.
\end{lemma}


\subsection{High-dimensional consistency of (AR)GES}\label{subsection: high-dimensional consistency of ARGES}


First, we define a collection of oracle versions of (AR)GES using the oracle scoring criterion $\mathcal{S}_{\lambda_n}^*$ (Definition \ref{definition: oracle score}). Every move in the forward or backward phase of (AR)GES from a CPDAG $\mathcal{C}_\text{current}$ to $\mathcal{C}_\text{new}$ corresponds to an edge addition or edge deletion in a DAG in the Markov equivalence class of $\mathcal{C}_{\text{current}}$, and thus corresponds to a partial correlation, by Lemma \ref{lemma: oracle score difference}. We denote the partial correlation associated with a move from $\mathcal{C}_\text{current}$ to $\mathcal{C}_\text{new}$ by $\rho(\mathcal{C}_\text{current},\mathcal{C}_\text{new})$.


At every step in the forward phase of (AR)GES, an optimal CPDAG is chosen among a set of possible choices for the next step (see line \ref{line: optimal choice 1} of Algorithm \ref{algorithm: the forward phase of ARGES}). These optimal choices in the forward phase are not crucial for consistency of (AR)GES \citep{ChickeringMeek02}. Thus, we define below oracle versions of (AR)GES that allow sub-optimal choices for edge additions. One of our assumptions (see (A5) below) will be based on this definition.

\begin{definition}\label{definition: oracle versions of ARGES} ($\delta$-optimal oracle version of (AR)GES)
Let $\delta \in [0,1]$. A $\delta$-optimal oracle version of (AR)GES with scoring criterion $\mathcal S^*_{\lambda}$ consists of two phases: a $\delta$-optimal oracle forward phase and an oracle backward phase. A $\delta$-optimal oracle forward phase of (AR)GES is Algorithm \ref{algorithm: the forward phase of ARGES} based on $\mathcal{I}_{0}$ (for ARGES-CIG) or $\skeleton(\mathcal{C}_{0})$ (for ARGES-skeleton) or the complete undirected graph (for GES), using the oracle scoring criterion $\mathcal{S}_{\lambda}^*$, where at each step with $\mathfrak{C} \neq \emptyset$, $\mathcal{C}_\text{new}$ is chosen to be any CPDAG in $\mathfrak{C}$ (see line \ref{line: optimal choice 1} of Algorithm \ref{algorithm: the forward phase of ARGES}) such that
\vspace{-0.1in}
$$ |\rho(\mathcal{C}_\text{current},\mathcal{C}_\text{new}) |  \geq \underset{\mathcal{C} \in \mathfrak{C}}{\max} |\rho(\mathcal{C}_\text{current},\mathcal{C}) | - \delta. \vspace{-0.05in}$$
An oracle backward phase of (AR)GES equals the backward phase of GES (Algorithm \ref{A2-algorithm: the backward phase of GES} of the supplementary material), 
using the oracle scoring criterion $\mathcal{S}_{\lambda}^*$. At each step with $\mathfrak{C} \neq \emptyset$, if there are several CPDAGs with the same optimal score, then one of these is chosen arbitrarily as $\mathcal{C}_\text{new}$ (see line \ref{A2-line: optimal choice 2} of Algorithm \ref{A2-algorithm: the backward phase of GES} of the supplementary material). 


\end{definition}

\begin{theorem}\label{theorem: soundness of oracle versions} (Soundness)
Assume that the distribution of $\mathbf{X}$ is multivariate Gaussian and DAG-perfect. Let $\delta \in [0,1]$. Let $m$ be such that the maximum degree in the output of the forward phase of every $\delta$-optimal oracle version of (AR)GES with scoring criterion $\mathcal{S}_{\lambda}^*$  is bounded by $m$ for all $\lambda \geq 0$. If $\lambda < - \frac{1}{2}\log(1 - \rho_{ij|S}^2)$ for all $i,j \in \{1,\ldots, p\}$ and $S \subseteq  \{1,\ldots, p\}\setminus\{i,j\}$ such that $|S| \leq m$ and $\rho_{ij|S} \neq 0$, then the outputs of all $\delta$-optimal oracle versions of (AR)GES with scoring criterion $\mathcal{S}_{\lambda}^*$ are identical and equal to $\mathcal{C}_{0}$.
\end{theorem}

Note that the edge additions in the forward phase of a $\delta$-optimal oracle version of (AR)GES are only slightly sub-optimal for small values of $\delta$. In fact, we let $\delta_n$ tend to zero as $n \rightarrow \infty$ in assumption (A5) below.  For $\delta = 0$, we refer to the forward phase of the $\delta$-optimal oracle version of (AR)GES as the oracle forward phase of (AR)GES.\\

We now consider an asymptotic scenario where the number of variables in $\mathbf{X}$ and the distribution of $\mathbf{X}$ are allowed to change with the sample size $n$. Thus, let $\{\mathbf{X}_n \}$ be a sequence of random vectors such that the distribution of each $\mathbf{X}_n$ is multivariate Gaussian and DAG-perfect. Further, we slightly modify Notation \ref{notation: true graphs} as follows.

\begin{notation}\label{notation: high-dim} We reserve the notation $F_n$, $\mathcal{G}_{n0}$, $\mathcal{C}_{n0}$, $\mathcal{I}_{n0}$ and $\mathcal{D}_n$ for the following: $F_n$ denotes the distribution of $\mathbf{X}_n = (X_{n1},\ldots,X_{np_n })^T$, $\mathcal{G}_{n0}$ denotes a perfect map of $F_n$,  $\mathcal{C}_{n0} = \CPDAG(\mathcal{G}_{n0})$ is the corresponding CPDAG,
  $\mathcal{I}_{n0}$ is the CIG of $F_n$,
  and $\mathcal D_n$ denotes the data, consisting of $n$ i.i.d.\ observations from $F_n$.
\end{notation}

We make the following assumptions to prove high-dimensional consistency of (AR)GES. 

\begin{description}
   \item[(A1)] (Gaussianity) The distribution of $\mathbf{X}_n$ is multivariate Gaussian and DAG-perfect for all $n$.
   
   \item[(A2)] (high-dimensional setting) $p_n = \mathcal{O}(n^{a})$ for some $0 \leq a < \infty$.
   
   \item[(A3)]\label{assumption 3} (sparsity condition) Let $q_n = \max_{1 \leq i \leq p_n} |\mathbf{Adj}_{\mathcal{C}_{n0}}(X_{ni})|$ be the maximum degree in $\mathcal{C}_{n0}$.  Then $q_n = \mathcal{O}(n^{1-b_1})$ for some $0 < b_1 \leq 1$.
     \item[(A4)] (consistent estimators of the CIG or the CPDAG-skeleton) There exists a sequence of estimators $\hat{\mathcal{I}}_n$ (for ARGES-CIG) or a sequence of estimators $\hat{\mathcal{U}}_n$ (for ARGES-skeleton) such that \vspace{-0.05in}$$\underset{n \to \infty}{\lim}~\Prob (\hat{\mathcal{I}}_n = \mathcal{I}_{n0}) = 1~\text{or}~ \underset{n \to \infty}{\lim}~\Prob (\hat{\mathcal{U}}_n = \skeleton(\mathcal{C}_{n0})) = 1.$$ 

	
	\item[(A5)] (bounds on the growth of oracle versions) The maximum degree in the output of the forward phase of every $\delta_n$-optimal oracle version of (AR)GES with scoring criterion $\mathcal S^*_{\lambda_n}$ is bounded by $K_nq_n$, for all $\lambda_n\ge 0$ and some sequences $\delta_n \rightarrow 0$ and $\{K_n\}$ such that $\delta_n^{-1} = \mathcal{O}(n^{d_1})$ and $K_n =  \mathcal{O}(n^{b_1-b_2})$ for some constants $b_2$ and $d_1$ satisfying $0 \leq d_1 < b_2/2 \leq 1/2$, where $q_n$ is given by (A3).
	
   \item[(A6)] (bounds on partial correlations) The partial correlations $\rho_{nij | S}$ between $X_{ni}$ and $X_{nj}$ given $\{X_{nr} : r \in S\}$ satisfy the following upper and lower bounds for all $n$, uniformly over $i,j \in \{1,\ldots,p_n \}$ and $S \subseteq \{1,\ldots,p_n\}\setminus \{ i,j\}$ such that $|S|\leq K_nq_n$ (where $K_n$ and $q_n$ are from (A3) and (A5)):\vspace{-0.05in}
   \begin{align*}
    \mathop{\sup}_{i\neq j, S} |\rho_{nij | S}| \leq M < 1,  ~~ \text{and} ~~  \mathop{\inf}_{i, j, S} \{|\rho_{nij | S}| : \rho_{nij | S} \neq 0 \} \geq c_n,\vspace{-0.05in}
   \end{align*}
   for some sequence $c_n \to 0$ such that $c_n^{-1} = \mathcal{O}(n^{d_2})$ for some constant $d_2$ satisfying $0 <d_2< b_2/2$, where $b_2$ is given by (A5).
%
\end{description}

\cite{KalischBuehlmann07a} proved high dimensional consistency of the PC algorithm assuming (A1), (A2), (A3) and a slightly weaker version of (A6). More precisely, the authors assumed (A6) with $d_2 <b_1/2$. The most criticized assumption among these four assumptions is probably (A6) \citep{UhlerEtAl13}, which is also known as the strong faithfulness condition. Interestingly, \cite{VandeGeerBuhlmann13} proved high-dimensional consistency of the global optimum of the $\ell_0$-penalized likelihood score function without assuming strong faithfulness, but assuming a permutation beta-min condition. Unfortunately, we cannot guarantee that a greedy search method like (AR)GES can always find the global optimum without (A6) and thus we cannot substitute the strong faithfulness assumption here. 

We make two additional assumptions compared to \cite{KalischBuehlmann07a}, namely (A4) and (A5). Assumption (A4) is natural and
it is not a strong assumption, since there are various estimation methods for the CIG or the CPDAG-skeleton that are consistent in sparse high-dimensional settings (e.g., \cite{MeinshausenBuehlmann06, BanerjeeEtAl08, FriedmanEtAl08, RavikumarEtAl08, RavikumarEtAl11, CaiEtAl12, KalischBuehlmann07a, HaEtAl15}). 
 We will discuss (A5) in Section \ref{subsection: assumption A5} below.





\begin{theorem}\label{theorem: high-dimensional consistency of ARGES}
Assume (A1) - (A6). Let $\hat{\mathcal C}_n$, $\breve{\mathcal C}_n$ and $\tilde{\mathcal C}_n$ be the outputs of ARGES-CIG based on $\hat{\mathcal{I}}_n$, ARGES-skeleton based on $\hat{\mathcal{U}}_n$ and GES respectively, with the scoring criterion $\mathcal{S}_{\lambda_n}$. Then there exists a sequence $\lambda_n \rightarrow 0$ such that \vspace{-0.1in}$$\underset{n \rightarrow \infty}{\lim}~\Prob(\hat{\mathcal{C}}_n = \mathcal{C}_{n0}) = \underset{n \rightarrow \infty}{\lim}~\Prob(\breve{\mathcal{C}}_n = \mathcal{C}_{n0})  = \underset{n \rightarrow \infty}{\lim}~\Prob(\tilde{\mathcal{C}}_n = \mathcal{C}_{n0}) = 1.$$
\end{theorem}
%
We choose $\lambda_n = \frac{1}{9}\log(1 - c_n^2)$ to prove the above theorems (see Section \ref{A2-subsection: high-dimensional consistency proof} of the supplementary material), where $c_n$ is as in (A6).
However, similar arguments hold for any sequence $\{\lambda_n\}$ satisfying $\lambda_n < \frac{1}{8}\log(1 - c_n^2)$ and $\frac{K_nq_n\log(p_n)}{n\lambda_n}\rightarrow 0$, where $p_n$, $q_n$ and $K_n$ are as in (A2), (A3) and (A5) respectively. The penalty parameter of the BIC criterion (see Remark \ref{remark: BIC}), which is known to be weak for sparse high-dimensional settings (see, for example, \cite{FoygelDrton10}), does not satisfy $\frac{K_nq_n\log(p_n)}{n\lambda_n}\rightarrow 0$, except for some very restricted choices of $\{p_n\}$, $\{q_n\}$ and $\{K_n\}$. Further, we note that \cite{VandeGeerBuhlmann13} proved their high-dimensional consistency result with $\lambda_n = \mathcal{O}(\sqrt{\log(p_n)/n})$. This choice of $\{\lambda_n\}$ satisfies $\frac{K_nq_n\log(p_n)}{n\lambda_n}\rightarrow 0$ under (A1) - (A6), when the constant $b_2$ in (A5) is greater than $1/2$.




\subsection{Discussion on assumption (A5)}\label{subsection: assumption A5}

Note that for every $\delta_n$-optimal forward phase of (AR)GES with scoring criterion $\mathcal{S}_{\lambda_n}^*$, the maximum degree in the output increases or remains unchanged as $\lambda_n$ decreases. Thus, without loss of generality, we fix the scoring criterion for all $\delta_n$-optimal versions of (AR)GES to be $\mathcal{S}_{0}^*$ in (A5) and in the remainder of this subsection.

We first consider (A5) for $\{\delta_n\} = 0$. Then (A3) and (A5) together imply that the output of the oracle forward phase of (AR)GES is bounded by $K_nq_n = \mathcal{O}(n^{1-b_2})$, where $K_n$ is an upper bound on the ratio of the maximum degrees in the output of the oracle forward phase of (AR)GES and in $\mathcal{C}_{n0}$. It follows from the soundness of the oracle version of (AR)GES that $K_n\ge 1$ and hence $b_2 \le b_1$, where $b_1$ is given by (A3). Therefore, (A5) roughly states that the oracle forward phase of (AR)GES does not add ``too many" edges that do not belong to the CPDAG-skeleton. In Section \ref{A2-section: verifying (A5)}
of the supplementary material, we empirically verify this in certain sparse high-dimensional settings and we find that in all but very few cases the maximum degree of the outputs of the oracle forward phase of (AR)GES are reasonably close to the maximum degree of the true CPDAG-skeleton. Further, note that (A5) is slightly different for GES, ARGES-CIG and ARGES-skeleton (see Figure \ref{A2-fig: boxplots of max degree} of the supplementary material).

The intuition for having a $\delta_n$-optimal oracle version in (A5), with $\delta_n >0$ is the following: If
$$
  |\rho(\mathcal{C}_\text{current},\mathcal{C}_\text{new}) |  \geq \underset{\mathcal{C} \in \mathfrak{C}}{\max} |\rho(\mathcal{C}_\text{current},\mathcal{C}) | - \delta_n
$$
(see Definition \ref{definition: oracle versions of ARGES}), then the move from $\mathcal{C}_\text{current}$ to $\mathcal{C}_\text{new}$ is so close to the optimal move that it cannot be identified as sub-optimal in the sample version.\\

We note that the output of a $\delta$-optimal oracle forward phase of (AR)GES depends not only on the structure $\mathcal{C}_{0}$ but also on the absolute values of the nonzero partial correlations (even for $\delta =0$). The latter makes it very difficult to characterize a family of distributions for which (A5) holds. However, we provide two strong structural conditions under which (A5) holds. The conditions are based on the following two results.

\begin{lemma}\label{lemma: disconnected nodes}
Assume that the distribution of $\mathbf{X}$ is multivariate Gaussian and DAG-perfect. Let $X_{i}$ and $X_{j}$ be two nodes in $\mathcal C_{0}$. If there is no path between $X_{i}$ and $X_{j}$ in $\mathcal C_{0}$, then the output of a $\delta$-optimal oracle forward phase of (AR)GES does not contain an edge between $X_{i}$ and $X_{j}$.
\end{lemma}

\begin{theorem}\label{theorem: forest}
Assume that the distribution of $\mathbf{X}$ is multivariate Gaussian and DAG-perfect. If $\skeleton(\mathcal{C}_{0})$ is a forest (i.e.\ contains no cycle), then the output of the oracle forward phase of (AR)GES equals $\mathcal{C}_{0}$.
\end{theorem}


The proof of Lemma \ref{lemma: disconnected nodes} is trivial, since if $X_{i}$ and $X_{j}$ are not connected by a path in $\mathcal{C}_{0}$, then $\rho_{ij|S} = 0$ for all $S \subseteq\{1,\ldots,p\}\setminus \{i,j\}$. Further, note that Lemma \ref{lemma: disconnected nodes} implies that it is sufficient to prove Theorem \ref{theorem: forest} when $\skeleton(\mathcal{C}_{0})$ is an undirected tree (i.e.\ a connected graph containing no cycle) or equivalently, when $\mathcal{G}_{0}$ is a polytree (i.e.\ a directed graph whose skeleton is an undirected tree).


We note that Theorem \ref{theorem: forest} shows a connection between the oracle forward phase of GES and the Chow-Liu algorithm \citep{ChowLiu68} for multivariate Gaussian distributions. The Chow-Liu algorithm is a greedy forward search procedure for learning optimal undirected trees, based on mutual information between pairs of variables (which equals $-\frac{1}{2}\log(1-\rho_{ij}^2)$ for a multivariate Gaussian distribution). Theorem \ref{theorem: forest} shows that the oracle forward phase of GES is a greedy forward search procedure for learning optimal polytrees, based on conditional mutual information $-\frac{1}{2}\log(1-\rho_{ij|S}^2)$. To our knowledge, this connection between GES and the Chow-Liu algorithm cannot be found in the literature on the Chow-Liu algorithm and extensions thereof for learning polytrees \citep{RebanePearl87, HuetedeCampos93, deCampos98, OuerdEtAl04}.

We are now ready to state our sufficient conditions for assumption (A5) as the following immediate corollaries of Lemma \ref{lemma: disconnected nodes} and Theorem \ref{theorem: forest}.

\begin{corollary}\label{corollary: disconnected nodes}
Assume (A1). If the number of nodes in each connected component of $\mathcal{C}_{n0}$ is bounded by $q_n'$ satisfying $q_n' = \mathcal{O}(n^{a'})$ for some $0 \leq a' <1$, then the maximum degree in the output of a $\delta_n$-optimal oracle forward phase of (AR)GES is bounded by $q_n'$ for all $\delta_n \geq 0$.
\end{corollary}

\begin{corollary}\label{corollary: forest}
Assume (A1) and (A3). If $\skeleton(\mathcal{C}_{n0})$ is a forest, then the maximum degree in the output of the oracle forward phase of (AR)GES equals $q_n$, where $q_n$ is given by (A3).
\end{corollary}

In order to extend Theorem \ref{theorem: forest} or Corollary \ref{corollary: forest} to a $\delta_n$-optimal version of (AR)GES, we need the following additional assumption.

\begin{description}
\item[(A7)] (bounds on the gaps between marginal correlations and first order partial correlations) \vspace{-0.1in}$$\underset{(i, j, k) \in T}{\inf}~  \{\left||\rho_{nij | k}| - |\rho_{nij}|\right| : |\rho_{nij | k}| \neq |\rho_{nij}|\}  \geq c_n',$$ where $T = \{(i,j,k) : \text{$(X_{ni},X_{nj},X_{nk})$ is an unshielded triple in $\mathcal{C}_{n0}$} \}$ and $c_n' \to 0$ such that $c_n^{\prime -1} = \mathcal{O}(n^{d_3})$ for some $0\leq d_3 < b_1/2$ where $b_1$ is given by (A3).
\end{description}

To understand (A7), assume that  $\skeleton(\mathcal{C}_{n0})$ is a forest containing an unshielded triple $(X_{ni},X_{nj},X_{nk})$ such that $\pi_{\mathcal{C}_{n0}}(X_{ni},X_{nj},X_{nk})$ is either $X_{ni} \to X_{nj} \leftarrow X_{nk}$ or $X_{ni} - X_{nj} - X_{nk}$. Let $\mathcal{C}_n$ be the CPDAG representing the Markov equivalence class of the DAG $\mathcal{G}_n$ obtained by deleting the edge between $X_{ni}$ and $X_{nj}$ in $\mathcal{G}_{n0}$. Note that $\rho_{nij | k}$ corresponds to the improvement in the oracle score for adding the edge $X_{ni} \to X_{nj}$ in a DAG in the Markov equivalence class of $\mathcal{C}_n$ while creating a new v-structure $X_{ni} \to X_{nj} \leftarrow X_{nk}$, whereas $\rho_{nij}$ corresponds to the improvement in the oracle score for adding the edge $X_{ni} \to X_{nj}$ in a DAG in the Markov equivalence class of $\mathcal{C}_n$ without creating the v-structure $X_{ni} \to X_{nj} \leftarrow X_{nk}$. Thus (A7) ensures that for each $\delta_n< c_n'$, a $\delta_n$-optimal oracle forward phase of (AR)GES would move from $\mathcal{C}_n$ to $\mathcal{C}_{n0}$ by correctly deciding if a new v-structure should be created.

\begin{theorem}\label{theorem: forest 2}
Assume (A1), (A3), (A6) with $K_n = 1$ and (A7). If $\skeleton(\mathcal{C}_{n0})$ is a forest, then the output of a $\delta_n$-optimal oracle forward phase of (AR)GES equals $\mathcal{C}_{n0}$ for all $\delta_n < \min((1-M)c_n,c_n')$, where $M$ and $c_n$ are given by (A6) and $c_n'$ is given by (A7). Hence, (A5) holds with $b_2 =b_1$ and $d_1 = \max(d_2,d_3)$, where $b_1$ is given by (A3).
\end{theorem}

\begin{remark}
An obvious extension of Corollary \ref{corollary: disconnected nodes} and Theorem \ref{theorem: forest 2} is the following: if the assumption of Corollary \ref{corollary: disconnected nodes} holds for the connected components of $\skeleton(\mathcal{C}_{n0})$ that are not trees and the assumptions of Theorem \ref{theorem: forest 2} hold for the connected components of $\skeleton(\mathcal{C}_{n0})$ that are trees, then (A5) holds with $b_2=a'$ and $d_1 = \max(d_2,d_3)$.
\end{remark}







\section{High-dimensional consistency of GES and ARGES for linear structural equation models}\label{section: high dimensional consistency for LSEM}


In this section, we present a slightly modified version of the result from the previous section: we prove high-dimensional consistency of (AR)GES for linear structural equation models with sub-Gaussian error variables.

\begin{definition}\label{definition: linear SEM} 
   Let $\mathcal{G}_0 = (\mathbf{X},E)$ be a DAG and let $B_{\mathcal{G}_0}$ be a $p \times p$ matrix such that $(B_{\mathcal{G}_0})_{ji} \neq 0$ if and only if $X_i \in \mathbf{Pa}_{\mathcal{G}_0}(X_j)$. Let $\boldsymbol{\epsilon} = (\epsilon_1,\ldots,\epsilon_p)^{T}$ be a random vector of jointly independent error variables. Then $\mathbf{X} = (X_1,\dots,X_p)^T$ is said to be generated from a linear structural equation model (linear SEM) characterized by $ (B_{\mathcal{G}_0},\boldsymbol{\epsilon})$ if $\mathbf{X} = B_{\mathcal{G}_0} \mathbf{X} + \boldsymbol{\epsilon}.$
\end{definition}

We assume the distribution of $\mathbf{X}$ is faithful to $\mathcal{G}_0$, implying that $\mathcal{G}_0$ is a perfect map of the distribution of $\mathbf{X}$. We refer to this as $\mathbf{X}$ is generated from a DAG-perfect linear SEM. 

Note that if $\epsilon_1,\ldots,\epsilon_p$ are Gaussian random variables, then the joint distribution of $\mathbf{X}$ is multivariate Gaussian with covariance matrix $\Sigma_{0} = (\mathrm{I} - B_{\mathcal{G}_0})^{-1}\Cov(\boldsymbol{\epsilon})(\mathrm{I} - B_{\mathcal{G}_0})^{-T}$. In this case, for any $i\neq j$ and $S \subseteq \{1,\ldots,p\} \setminus \{i,j\}$, we have
\vspace{-0.05in}
\begin{align}\label{eq: d-sep and partial correlation}
X_i \dsep_{\mathcal{G}_0} X_j \mid \{X_r : r \in S\} \iff \rho_{ij|S} = 0.
\end{align}
Since $\rho_{ij|S} = 0$ depends on the distribution of $\mathbf{X}$ via $\Sigma_0$, \eqref{eq: d-sep and partial correlation} holds regardless of the distribution of $\boldsymbol{\epsilon}$ (see \citet{SpirtesEtAl98}). Therefore, if $\mathbf{X}$ is generated from a DAG-perfect linear SEM, then $X_i \indep X_j \mid \{X_r : r \in S\}$ if and only if $\rho_{ij|S} = 0$. Consequently, Lemma \ref{lemma: oracle score difference} and Theorem \ref{theorem: soundness of oracle versions} imply the soundness of the (AR)GES algorithm with the $\ell_0$-penalized Gaussian log-likelihood oracle score $\mathcal{S}_{\lambda}^*(\cdot,\tilde{F})$ (Definition \ref{definition: oracle score}) for a DAG-perfect linear SEM characterized by $ (B_{\mathcal{G}_0},\boldsymbol{\epsilon})$, where $\tilde{F}$ is the distribution of a Gaussian random vector with covariance matrix $\Sigma_{0} = (\mathrm{I} - B_{\mathcal{G}_0})^{-1}\Cov(\boldsymbol{\epsilon})(\mathrm{I} - B_{\mathcal{G}_0})^{-T}$. 

We use Notation \ref{notation: high-dim} to present the following high-dimensional consistency result for linear SEMs with sub-Gaussian error variables. We replace (A1) given in Section \ref{subsection: high-dimensional consistency of ARGES} by (A1*) below, and we make an additional assumption.
\begin{description}
\item[(A1*)] For each $n$, $\mathbf{X}_n$ is generated from a DAG-perfect linear SEM with sub-Gaussian error variables satisfying $||\epsilon_{ni}||_{\psi_2} \leq C_1 \sqrt{\Var(\epsilon_{ni})}$ for some absolute constant $C_1>0$, where $||\cdot||_{\psi_2}$ denotes the sub-Gaussian norm.
\item[(A8)] For any $(K_nq_n+2) \times (K_nq_n+2)$ principal submatrix $\Sigma_n$ of $\Sigma_{n0} = \Cov(\mathbf{X}_n)$, \vspace{-0.05in}$$C_2 \leq 1/||\Sigma_n^{-1}||_2 \leq ||\Sigma_n||_2 \leq C_3$$ for some absolute constants $C_2,C_3>0$, where $||\cdot||_2$ denotes the spectral norm and $q_n$ and $K_n$ are given by (A3) and (A5) respectively.
\end{description}

\begin{theorem}\label{theorem: LSEM high-dimensional consistency}
Assume (A1*), (A8) and (A2) - (A6) given in Section \ref{subsection: high-dimensional consistency of ARGES}.
Then there exists a sequence $\lambda_n \rightarrow 0$ such that $$\underset{n\to\infty}{\lim}\Prob(\hat{\mathcal{C}}_n = \mathcal{C}_{n0}) = \underset{n\to\infty}{\lim}\Prob(\breve{\mathcal{C}}_n = \mathcal{C}_{n0}) = \underset{n\to\infty}{\lim}\Prob(\tilde{\mathcal{C}}_n = \mathcal{C}_{n0}) =1,$$
where $\hat{\mathcal C}_n$, $\breve{\mathcal C}_n$ and $\tilde{\mathcal C}_n$ are the outputs of ARGES-CIG and ARGES-skeleton and GES respectively, with the $\ell_0$-penalized Gaussian log-likelihood scoring criterion $\hat{\mathcal{S}}_{\lambda_n}$ (see Definition \ref{definition: penalized likelihood score}).
\end{theorem}
\begin{remark}
As in \cite{HarrisDrton13}, we can allow an ultra high-dimensional setting by replacing (A2) by $p_n = \mathcal{O}(\exp(n^{a}))$ (for some $0\leq a <1$) and allow an $\mathcal{O}(n^{f})$ (for some $0\leq f <1/4$) growth rate of $||\Sigma_n^{-1}||_2$ for each $(K_nq_n+2) \times (K_nq_n+2)$ submatrix $\Sigma_n$ of $\Sigma_{n0}$ under the additional restriction that $a+ 4f < b_2 - 2\max(d_1,d_2) $.
\end{remark}

\begin{remark}
We note that the results of Section \ref{subsection: assumption A5} continue to hold if we replace the multivariate Gaussian assumption by the assumption that $\mathbf{X}_n$ is generated from a linear SEM with arbitrary error variables.
\end{remark}

\begin{remark}
We note that a similar high-dimensional consistency result holds for the PC algorithm when the sample partial correlations are used for testing conditional independence.
\end{remark}

\section{High-dimensional consistency of GES and ARGES  in the nonparanormal setting}\label{section: high dimensional consistency in the NPN setting}

In this section, we prove high-dimensional consistency under the assumption that each $\mathbf{X}_n$ has a nonparanormal distribution.

\begin{definition}\label{definition: nonparanormal} \citep{HarrisDrton13}
Let $\Sigma_0 = (\rho_{ij})$ be a positive definite correlation matrix and let $\mathbf{Z} = (Z_1,\ldots,Z_p)^T \sim N(0, \Sigma_0)$.  For a collection of strictly increasing functions $\mathbf{g} = ( g_1,\ldots,g_{p})^T$, the random vector $\mathbf{X} = ( g_1(Z_1),\ldots,g_{p}(Z_p))^T$ has a nonparanormal distribution $\NPN(\mathbf{g} , \Sigma_0)$.
\end{definition}

A variant of the PC algorithm, called Rank PC, was shown to be consistent in high-dimensional settings under the nonparanormal distribution assumption \citep{HarrisDrton13}. First, we briefly discuss the Rank PC algorithm. Then we define a new scoring criterion motivated by the rank-based measures of correlations used in the rank PC algorithm. Finally, we prove high-dimensional consistency of (AR)GES with this scoring criterion in the nonparanormal setting.

Let $\Sigma_0$, $\mathbf Z$, $\mathbf g$, and $\mathbf{X}$ be as in Definition \ref{definition: nonparanormal}. First, note that since the marginal transformations $g_i$ are deterministic, for any $i \neq j$ and $S \subseteq \{1,\ldots, p \} \setminus \{i,j\}$
\vspace{-0.05in}
\begin{align}\label{eq: conditional independence in NPN}
X_i \indep X_j \mid \{X_r : r \in S\} \iff Z_i \indep Z_j \mid \{Z_r : r \in S\} \iff \rho_{ij|S} = 0,
\end{align}
where $\rho_{ij|S}$ is the partial correlation coefficient between $Z_i$ and $Z_j$ given $\{ Z_{r} : r \in S \}$. Next, note that since $g_i$ are strictly increasing functions, a sample rank correlation coefficient (Spearman's $\rho$ or Kendall's $\tau$) between $X_i$ and $X_j$ (denoted as $\hat{\rho}^{S}_{ij}$ or $\hat{\rho}^K_{ij}$) is identical to the corresponding rank correlation between $Z_i$ and $Z_j$. Further, \cite{LiuEtAl12} showed that $2\sin(\frac{\pi}{6} \hat{\rho}^{S}_{ij})$ and $\sin(\frac{\pi}{2} \hat{\rho}^K_{ij})$ are consistent estimators of $\rho_{ij}$.
In the remainder of this section, we generically denote a rank based estimator $2\sin(\frac{\pi}{6} \hat{\rho}^{S}_{ij})$ or $\sin(\frac{\pi}{2} \hat{\rho}^K_{ij})$ by $\hat{\rho}_{ij}$. We denote the corresponding estimator of the correlation matrix $\Sigma_0$ by $\hat{\Sigma}= (\hat{\rho}_{ij})$. Following \cite{HarrisDrton13}, we define rank based estimators of partial correlations $\rho_{ij|S}$ through the matrix inversion formula:
\vspace{-0.1in}
\begin{align}\label{eq: recursive formula for partial correlations}
\hat{\rho}_{ij|S} :=  - \frac{\hat{\Psi}_{12}^{-1}}{\sqrt{\hat{\Psi}_{11}^{-1} \hat{\Psi}_{22}^{-1} }},
\end{align}
where $\hat{\Psi}$ is the submatrix of $\hat{\Sigma}$ that corresponds to $X_{i}$, $X_{j}$ and $\{X_r : r \in S\}$ in this order, and $\hat{\Psi}_{st}^{-1}$ denotes the $(s,t)$th entry of $\hat{\Psi}^{-1}$.

Given the estimators $\hat{\rho}_{ij|S}$ defined above, the rank PC algorithm is just the PC algorithm based on the following conditional independence tests: reject the null hypothesis $X_i \indep X_j \mid \{X_r : r \in S\}$ if and only if $|\hat{\rho}_{ij|S}| > \nu$, where the critical value $\nu$ is chosen to be the same for each individual test, and it is a tuning parameter of the algorithm that controls sparsity of the output.

We now define the scoring criterion $\tilde{\mathcal{S}}_{\lambda_n}$ (Definition \ref{definition: NPN score}), which is motivated by Lemma \ref{lemma: score difference}. We show below that (AR)GES with $\tilde{\mathcal{S}}_{\lambda_n}$ is consistent in certain sparse high-dimensional settings with nonparanormal distributions.


\begin{definition}\label{definition: NPN score}
We define the scoring criterion $\tilde{\mathcal{S}}_{\lambda_n}$ by setting the score of an empty DAG to zero, and by defining the score difference between two DAGs that differ by exactly one edge as follows. Let $\mathcal{H}_n = (\mathbf{X}_n,E_n)$ be a DAG such that $X_{ni}\in \mathbf{Nd}_{\mathcal{H}_n}(X_{nk})\setminus \mathbf{Pa}_{\mathcal{H}_n}(X_{nk})$. Let $\mathcal{H}_n'=(\mathbf{X}_n,E_n \cup \{X_{ni} \to X_{nk}\})$. Then
\vspace{-0.05in}
\begin{align}\label{eq: NPN score difference}
\mathcal{\tilde{\mathcal{S}}}_{\lambda_n}(\mathcal{H}_n',\mathcal{D}_n) - \mathcal{\tilde{\mathcal{S}}}_{\lambda_n}(\mathcal{H}_n,\mathcal{D}_n) := \frac{1}{2}\log \left(1 - \hat{\rho}_{nik|\mathbf{Pa}_{\mathcal{H}_n}(k)}^2\right) + \lambda_n,
\end{align}
where $\mathbf{Pa}_{\mathcal{H}_n}(k) = \{j : X_{nj} \in \mathbf{Pa}_{\mathcal{H}_n}(X_{nk})\}$ and $\hat{\rho}_{nik|\mathbf{Pa}_{\mathcal{H}_n}(k)}$ is a rank correlation based estimate defined by \eqref{eq: recursive formula for partial correlations}. 

For a DAG $\mathcal{H}_n = (\mathbf{X}_n,E)$, the score $\mathcal{\tilde{\mathcal{S}}}_{\lambda_n}(\mathcal{H}_n,\mathcal{D}_n)$ can be obtained by summing up the score differences while sequentially adding directed edges from $E$ starting from the empty graph.
\end{definition}

\begin{lemma}\label{lemma: well-defined score}
$\mathcal{\tilde{\mathcal{S}}}_{\lambda_n}(\mathcal{H}_n,\mathcal{D}_n)$ is well-defined for all DAGs $\mathcal{H}_n$, i.e., the score does not depend on the order in which the directed edges are added to the empty graph. 

\end{lemma}

\begin{lemma}\label{lemma: NPN score equivalence}
$\mathcal{\tilde{\mathcal{S}}}_{\lambda_n}$ is score equivalent, i.e., $\mathcal{\tilde{\mathcal{S}}}_{\lambda_n}(\mathcal{H}_n',\mathcal{D}_n) = \mathcal{\tilde{\mathcal{S}}}_{\lambda_n}(\mathcal{H}_n,\mathcal{D}_n)$ for any two Markov equivalent DAGs $\mathcal{H}_n$ and $\mathcal{H}_n'$, and for all $\mathcal{D}_n$.
\end{lemma}


For a nonparanormal distribution $\NPN(\mathbf{g}_n , \Sigma_{n0})$, we define the oracle score $\mathcal{\tilde{\mathcal{S}}}_{\lambda_n}^*(\mathcal{H}_n,\Sigma_{n0})$ by replacing \eqref{eq: NPN score difference} in Definition \ref{definition: NPN score} with the following:
\begin{align*}
\mathcal{\tilde{\mathcal{S}}}_{\lambda_n}^*(\mathcal{H}_n',\Sigma_{n0}) - \mathcal{\tilde{\mathcal{S}}}_{\lambda_n}^*(\mathcal{H}_n,\Sigma_{n0}) := \frac{1}{2}\log \left(1 - \rho_{nik|\mathbf{Pa}_{\mathcal{H}_n}(k)}^2\right) + \lambda_n,
\end{align*}
where the partial correlations are given by $\Sigma_{n0}$. In fact, $\mathcal{\tilde{\mathcal{S}}}_{\lambda_n}^*(\cdot,\Sigma_{n0})$ is identical to the scoring criterion $\mathcal{\mathcal{S}}_{\lambda_n}^*(\cdot,\tilde{F}_n)$ (Definition \ref{definition: oracle score}), where $\tilde{F}_n$ is the distribution of $\mathbf{Z}_n \sim N(\mathbf{0},\Sigma_{n0})$.

\begin{theorem}\label{theorem: NPN high-dimensional consistency}
Assume that the distribution of $\mathbf{X}_n$ is $\NPN(\mathbf{g}_n , \Sigma_{n0})$ and DAG-perfect. Assume (A2) - (A6) given in Section \ref{subsection: high-dimensional consistency of ARGES}, with $\delta_n$-optimal oracle versions based on $\mathcal{\tilde{\mathcal{S}}}_{\lambda_n}^*$ (for (A5)), and 
with partial correlations based on $\Sigma_{n0}$ (for (A6)). 
Further, assume that the constants $b_1,b_2,d_1,d_2$ in (A3) - (A6) satisfy the following stronger restrictions: $1/2<b_2 \leq b_1\leq 1$ and $\max(d_1,d_2) <b_2 - 1/2$. Finally, assume that $1/||\Sigma_n^{-1}||_2$ is bounded below by an absolute constant $C_2>0$ for all $(K_nq_n+2) \times (K_nq_n+2)$ principal submatrices $\Sigma_n$ of $\Sigma_{n0}$, where $q_n$ and $K_n$ are given by (A3) and (A5) respectively. Then there exists a sequence $\lambda_n \rightarrow 0$ such that $$\underset{n\to\infty}{\lim}\Prob(\hat{\mathcal{C}}_n = \mathcal{C}_{n0}) = \underset{n\to\infty}{\lim}\Prob(\breve{\mathcal{C}}_n = \mathcal{C}_{n0}) = \underset{n\to\infty}{\lim}\Prob(\tilde{\mathcal{C}}_n = \mathcal{C}_{n0}) =1,$$
where $\hat{\mathcal C}_n$, $\breve{\mathcal C}_n$ and $\tilde{\mathcal C}_n$ are the outputs of ARGES-CIG, ARGES-skeleton and GES, respectively, with the scoring criterion $\tilde{\mathcal{S}}_{\lambda_n}$ given by Definition \ref{definition: NPN score}.
\end{theorem}

Our assumptions in Theorem \ref{theorem: NPN high-dimensional consistency} are similar to the assumptions of \cite{HarrisDrton13}\footnote{As in \cite{HarrisDrton13}, we can allow an ultra high-dimensional setting by replacing (A2) by $p_n = \mathcal{O}(\exp(n^{a}))$ (for some $0\leq a <1$) and allow an $\mathcal{O}(n^{f})$ (for some $0\leq f <1/4$) growth rate of $||\Sigma_n^{-1}||_2$ for each $(K_nq_n+2) \times (K_nq_n+2)$ submatrix $\Sigma_n$ of $\Sigma_{n0}$ if we replace the condition $\max(d_1,d_2) <b_2 - 1/2$ by $a+2(1-b_2)+2\max(d_1,d_2) + 4f < 1$.}, except that we additionally assume (A4) and (A5), where we require (A4) only for ARGES. Note that (A4) is not a strong assumption since there are high-dimensionally consistent estimators of the CIG or CPDAG-skeleton \citep{LiuEtAl12, HarrisDrton13}. In fact, \cite{LiuEtAl12} proposed to use rank based estimators of the correlation matrix as described above, and showed that high-dimensional consistency can be retained in the nonparanormal setting by plugging in such an estimated correlation matrix in many CIG estimation methods that are developed for multivariate Gaussian distributions (e.g. \cite{FriedmanEtAl08, CaiEtAl12}).






\section{Simulations}\label{section: high-dimensional simulations}

Having shown that (AR)GES has similar theoretical guarantees as the PC algorithm in high-dimensional settings, we now compare the finite sample performance and computational efficiency of (AR)GES and PC. In fact, we compare (AR)GES with the order independent version of PC \citep{ColomboMaathuis14}, but in the remaining of this section, we simply refer to it as PC. Additionally, we include another popular hybrid structure learning method, called Max-Min Hill-Climbing (MMHC) \citep{TsamardinosEtAl06}. MMHC first estimates the CPDAG-skeleton by applying the Max-Min Parents and Children (MMPC) algorithm \citep{TsamardinosEtAl06}, and then performs a hill-climbing DAG search on the space restricted to the estimated CPDAG-skeleton (see also Remark \ref{remark: inconsistency of some other hybrid algorithms}). We use the \texttt{R}-package \textbf{pcalg} \citep{KalischEtAl12} for (AR)GES and PC and use the \texttt{R}-package \textbf{bnlearn} \citep{Scutari10} for MMHC. In fact, we use a slightly modified version of (AR)GES that additionally includes a \emph{turning phase} \citep{Chickering02, HauserBuhlmann12} and an iteration over all three phases (see Section \ref{A2-section: turing phase} of the supplementary material for details).


\subsection{Simulation settings}\label{subsection: simulation settings}



For each of the four settings given in Table \ref{table: simulation settings}, we use the \texttt{R}-package \textbf{pcalg} \citep{KalischEtAl12} to simulate $r = 100$ random weighted DAGs $\{\mathcal{G}_n^{(1)},\ldots,\mathcal{G}_n^{(r)}\}$ with $p_n$ vertices and expected number of edges $e_n$, where each pair of nodes in a randomly generated DAG has the probability $e_n/\binom{p_n}{2}$ of being adjacent. The edge weights are drawn independently from a uniform distribution on $(-1,-0.1)\cup(0.1,1)$.


\begin{table}[ht]
\centering
\footnotesize
\caption{Simulation settings.}
\begin{tabular}{c|rrrr}
$n$ & 100 & 200 & 300 & 400 \\ \hline
  $p_n$ & 300 & 600 & 1200 & 2400 \\ \hline
  $e_n$ & 300 & 840 & 2100 & 4800
  \end{tabular}

\label{table: simulation settings}
\end{table}

Let $B_n^{(t)}$ denote the weight matrix of the weighted DAG $\mathcal{G}_n^{(t)}$, i.e.\ $(B_n^{(t)})_{ij} \neq 0$ if and only if the edge $X_j \to X_i$ is present in $\mathcal{G}_n^{(t)}$ and it then equals the corresponding edge weight. For $t = 1,\ldots, r$, the weight matrix $B_n^{(t)}$ and a random vector ${\boldsymbol\epsilon}_n^{(t)} = (\epsilon_{n1}^{(t)},\ldots,\epsilon_{np_n}^{(t)})^T$ define a distribution on $\mathbf{X}_n^{(t)} = (X_{1n}^{(t)},\ldots,X_{np_n}^{(t)})^T$ via the linear structural equation model $\mathbf{X}_n^{(t)} = B_n^{(t)}\mathbf{X}_n^{(t)} + {\boldsymbol\epsilon}_n^{(t)}$. We choose $\epsilon_{n1}^{(t)},\ldots,\epsilon_{np_n}^{(t)}$ to be zero mean Gaussian random variables with variances independently drawn from a $\mathrm{Uniform}[1,2]$ distribution. We aim to estimate $\CPDAG(\mathcal{G}_n^{(t)})$ from $n$ i.i.d.\ samples from the multivariate Gaussian distribution of $\mathbf{X}_n^{(t)}$.

\subsection{Estimation of the CIG and the CPDAG-skeleton for ARGES}\label{subsection: MB and MMPC}

We estimate the CIG for ARGES-CIG using neighborhood selection with the LASSO\footnote{We will later use adaptive LASSO \citep{Zou06} for this step (see Section \ref{subsection: dense graphs}).} of \cite{MeinshausenBuehlmann06}, where we use the implementation in the \texttt{R}-package \textbf{huge} \citep{ZhaoEtAl12}. Neighborhood selection involves a tuning parameter $\gamma_n$ that corresponds to the LASSO penalization, where larger values of $\gamma_n$ yield sparser estimated graphs. 
We choose $\gamma_n =$ 0.16, 0.14, 0.12 and 0.10 for $p_n=$ 300, 600, 1200 and 2400 respectively. In Section \ref{A2-section: additional simulations} of the supplementary material, we empirically investigate the influence of $\gamma_n$ on the performance of ARGES-CIG. We find that the performance of ARGES-CIG is not very sensitive to the choice of $\gamma_n$ in this simulation setting, as long as $\gamma_n$ is reasonably small (i.e., the estimated CIG is reasonably dense). In particular, estimation quality can be slightly improved by choosing a smaller $\gamma_n$, but with a loss of computational efficiency.

We estimate the CPDAG-skeleton for ARGES-skeleton using the MMPC algorithm of \cite{TsamardinosEtAl06}, where we use the implementation in the \texttt{R}-package \textbf{bnlearn} \citep{Scutari10}. MMPC involves a tuning parameter $\kappa_n$ that corresponds to the significance level of the conditional independence tests, where smaller values of $\kappa_n$ yield sparser estimated graphs. We choose $\kappa_n =$ 0.2, 0.15, 0.10 and 0.05 for $p_n=$ 300, 600, 1200 and 2400 respectively. We note that the MMPC algorithm is computationally expensive for large values of $\kappa_n$. 

\subsection{Results}\label{subsection: results}
As the scoring criterion for (AR)GES and MMHC, we use the $\ell_0$-penalized likelihood score (see Definition \ref{definition: penalized likelihood score}) with a number of choices for the penalty parameter $\lambda_n$. Similarly, we apply PC with a number of choices for its tuning parameter $\alpha_n$ (the significance level for conditional independence tests). Finally, we compare their estimation quality with averaged receiver operating characteristic (ROC) curves, where we average true positive rates and false positive rates for each value of the tuning parameters over $r$ iterations (cf.\ threshold averaging of ROC curves \citep{Fawcett06}). We do not apply GES and MMHC for the case $p_n = 2400$, since they are too slow to handle such large graphs. We additionally apply ARGES-CIG with the true CIG, and we call it ARGES-CIG*.



\begin{figure}[!t]
\centering
\includegraphics[width=\textwidth]{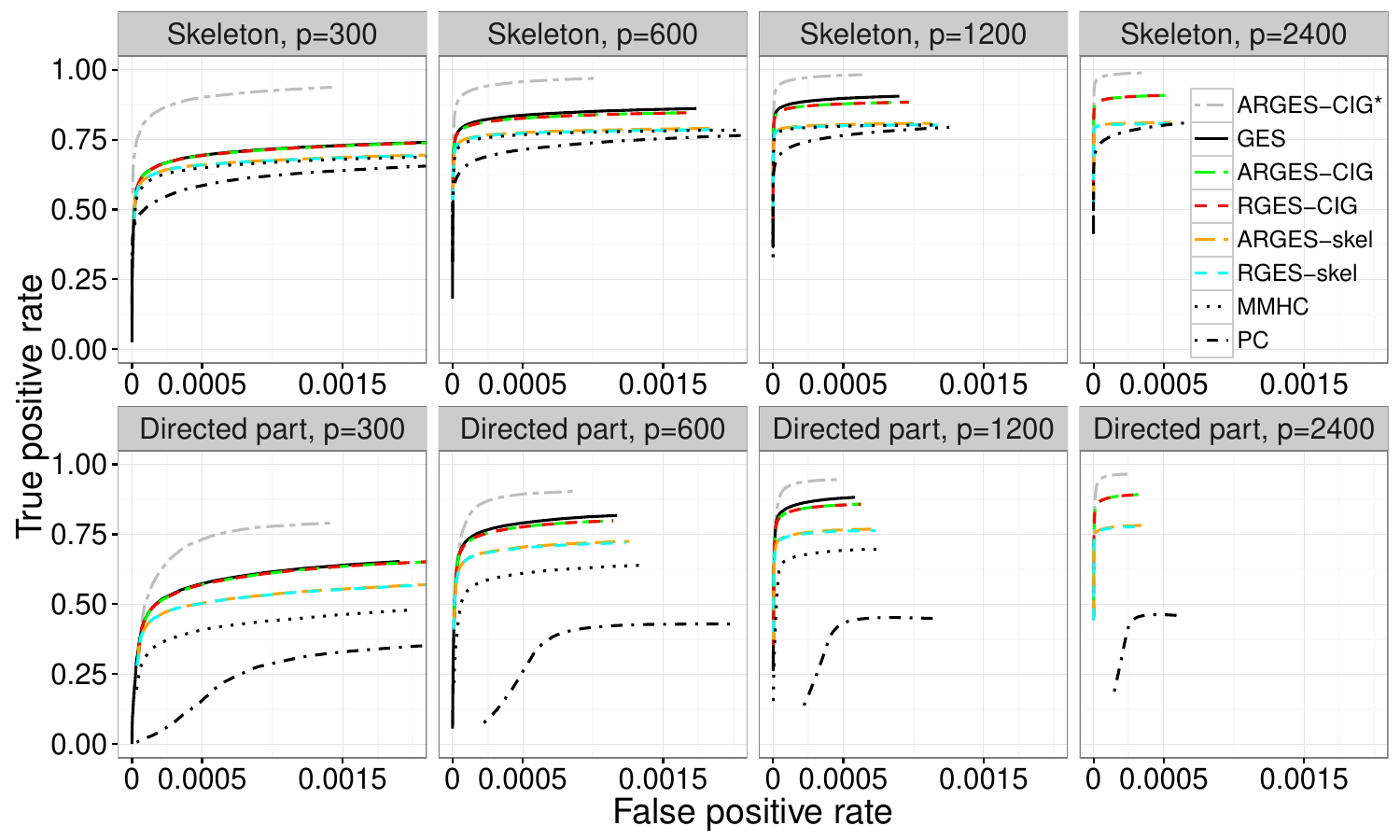}
\caption{Averaged ROC curves for estimating the skeletons (upper panel) and the directed parts (lower panel) of the underlying CPDAGs with ARGES-CIG*, ARGES-CIG, ARGES-skeleton, GES, MMHC and PC, for simulation settings given in Table \ref{table: simulation settings}.}
\label{fig: ROC curves}
\end{figure}

\begin{figure}[!t]
\centering
\includegraphics[width=0.75\textwidth]{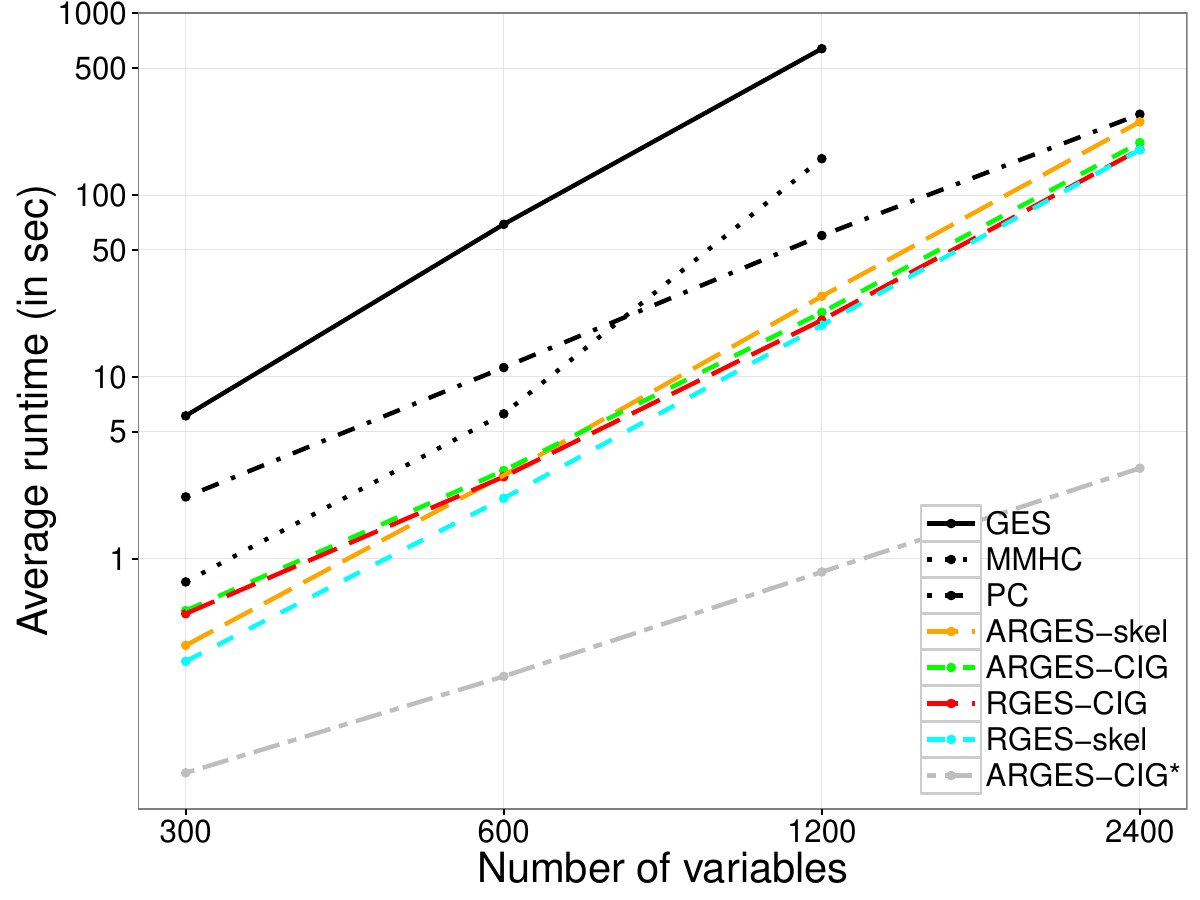}
\caption{Average runtimes (in seconds) for ARGES-CIG*, ARGES-CIG, ARGES-skeleton, GES, MMHC and PC, where the tuning parameters $\alpha_n$ and $\lambda_n$ are chosen to get roughly the right sparsity. 
The runtimes of (A)RGES do not include the CIG or the CPDAG-skeleton estimation part (see Table \ref{A2-table: CIG and MMPC outputs} of the supplementary material).}
\label{fig: runtime}
\end{figure}

In Figure \ref{fig: ROC curves}, we see that the averaged ROC curves get better as $n$ and $p_n$ increase together.
This supports the high-dimensional consistency theory of (AR)GES and PC. 
Based on Figure \ref{fig: ROC curves}, the performances of the algorithms can be summarized as follows:
\begin{align*}
\text{PC} < \text{MMHC} < \text{(A)RGES-skeleton} < \text{(A)RGES-CIG} \approx \text{GES} < \text{ARGES-CIG*},
\end{align*}
where $\text{A}<\text{B}$ represents that B outperformed A and $\text{A}\approx\text{B}$ represents that A and B performed similarly. 
Below, we list our main findings and possible explanations in detail.

\begin{enumerate}
\item ARGES-CIG* is the best performing method, but it is infeasible in practice as it requires knowledge of the true CIG.
\item GES is the next best performing method, closely followed by (A)RGES-CIG. However, the fact that ARGES-CIG* outperforms GES indicates the possibility that ARGES-CIG can outperform GES when combined with a better CIG estimation technique. We explore this in Section \ref{subsection: dense graphs}.
\item The fact that the performances of ARGES and RGES are almost identical, indicates that the adaptive part of ARGES does not have a significant influence on the performance.
\item (A)RGES-CIG outperformed (A)RGES-skeleton because (i) the true positive rate of (A)RGES is approximately bounded by the estimated CIG or CPDAG-skeleton, and (ii) the true positive rates of the estimated CIGs are larger than that of the estimated CPDAG-skeletons (see Table \ref{A2-table: CIG and MMPC outputs} of the supplementary material).
\item Although the performance of MMHC is similar to ARGES-skeleton for estimating the CPDAG-skeleton, it is significantly worse for estimating the directed part of the CPDAG. We suspect that this is due to some arbitrary choices of edge orientations, made in hill-climbing DAG search to resolve its non-uniqueness (see Section \ref{A2-subsection: limit outputs} of the supplementary material).
\item The constraint-based PC is the worst performing method in terms of estimation quality in these simulations. 
\end{enumerate}

Figure \ref{fig: runtime} shows the average runtimes of the algorithms\footnote{We use implementations of the algorithms from several \texttt{R}-packages. Hence, Figure \ref{fig: runtime} does not represent the computational efficiency of the algorithms, but rather represents a comparison of their currently available implementations in \texttt{R}.}. We see that GES does not scale well to large graphs and that the runtimes of (AR)GES and PC are somewhat similar. The runtimes of (A)RGES do not include the CIG or the CPDAG-skeleton estimation (see Table \ref{A2-table: CIG and MMPC outputs} of the supplementary material). The average runtimes of MMHC are much worse than the runtimes of (A)RGES-skeleton, although they are based on the same estimated CPDAG-skeleton. 

\subsection{Simulation with decreasing sparsity level}\label{subsection: dense graphs}


We have seen that in very sparse high-dimensional settings where GES performed very well, ARGES-CIG based on neighborhood selection could not outperform GES in terms of estimation.  One would expect, however, that using a good restricted search space can also be beneficial for the estimation performance. We investigate this in the following simulations, where we increased the adversity of the the problem by decreasing the sparsity level, and we paid more attention to the estimation of the CIG by including adaptive LASSO \citep{Zou06}. Moreover, to make the method applicable in practice, we no longer choose some pre-specified $\gamma_n$'s (cf.\ Section \ref{subsection: MB and MMPC}), but choose it via cross-validation.


\begin{table}[!h]
\centering
\footnotesize
\caption{Simulation settings with fixed number of variables and decreasing sparsity.}
\begin{tabular}{c|rrrr}
$n$ & 50 & 100 & 150 & 200 \\ \hline
  $p_n$ & 100 & 100 & 100 & 100 \\ \hline
  $e_n$ & 100 & 200 & 300 & 400
  \end{tabular}
\label{table: simulation settings 2}
\end{table}



\begin{figure}[!ht]
\centering
\includegraphics[width=\textwidth]{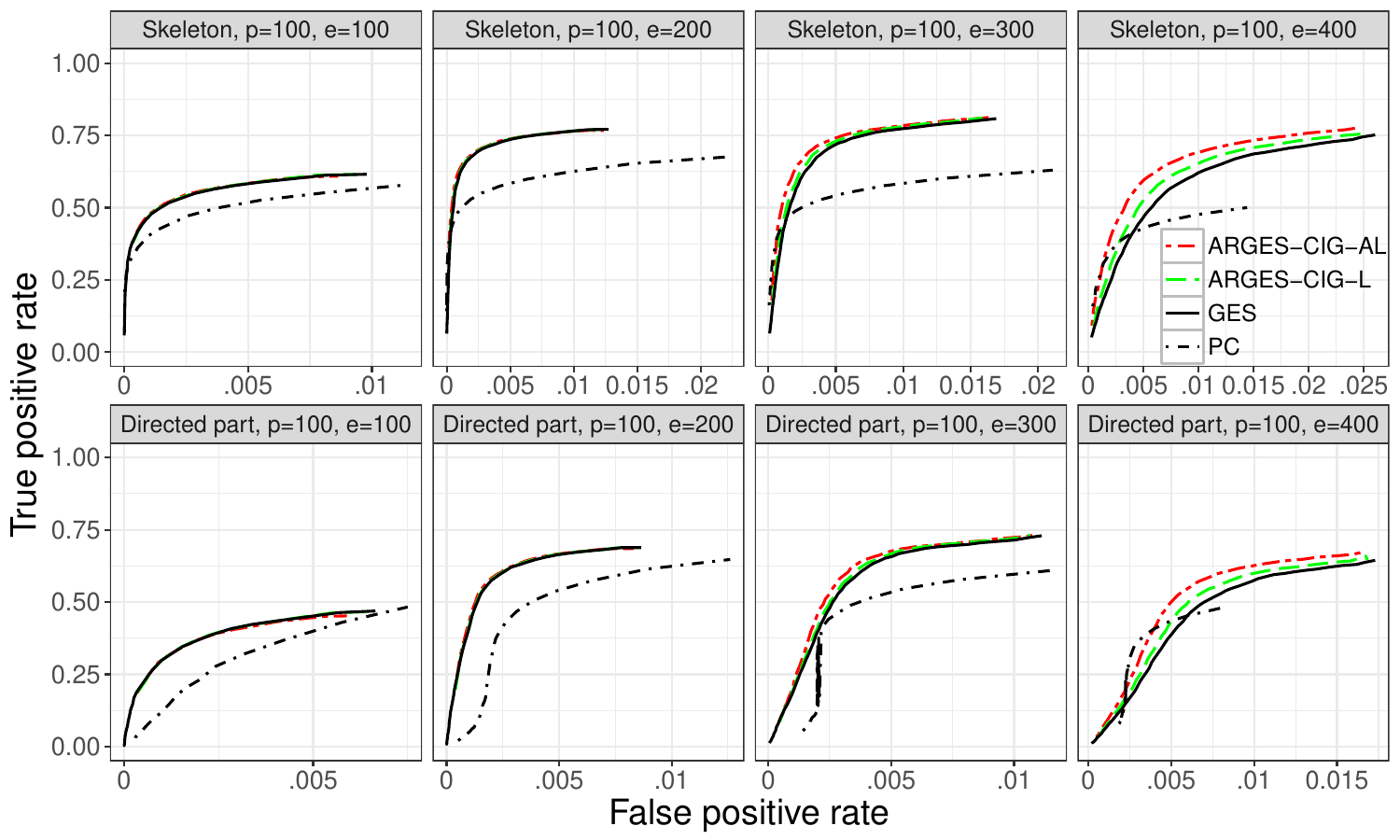}
\caption{Averaged ROC curves for estimating the skeletons (upper panel) and the directed parts (lower panel) of the underlying CPDAGs with ARGES-CIG based on adaptive LASSO (ARGES-CIG-AL), ARGES-CIG based on LASSO (ARGES-CIG-L), GES and PC, for simulation settings given in Table \ref{table: simulation settings 2}.}
\label{fig: ROC curves 2}
\end{figure}

We compare the following methods: ARGES-CIG based on LASSO, ARGES-CIG based on adaptive LASSO, GES and PC. We choose the tuning parameter for each LASSO\footnote{We use the implementation in the \texttt{R}-package \textbf{glmnet} \citep{FriedmanEtAl10}.} by 10-fold cross-validation, optimizing the prediction error (e.g., \citep{HastieEtAlbook09}). For each adaptive LASSO\footnote{We use the aforementioned implementation of LASSO after removing the variables with infinite weights and rescaling the other variables with the corresponding weights.}, we use the weights $w_i = 1/\hat{\beta}_i$ where $\hat{\beta}_i$ is the initial estimate of the $i$-th linear regression coefficient given by the corresponding LASSO regression with a 10-fold cross validation. After fixing the weights, we choose the tuning parameter in the second stage of the adaptive LASSO based on an additional 10-fold cross validation. ROC curves for estimating the skeleton of the CPDAG and the directed part of the CPDAG are obtained by varying $\lambda_n$ for all GES based methods, and by varying $\alpha_n$ for the PC algorithm.

Figure \ref{fig: ROC curves 2} shows that both versions of ARGES-CIG and GES perform equally well for the first three settings, while ARGES-CIG outperforms GES in the most adverse setting with $e_n=400$. Interestingly, although the estimated CIG based on the adaptive LASSO is always a subgraph of the estimated CIG based on the LASSO, the use of the adaptive LASSO for the CIG estimation enhances the performance gain as the sparsity decreases (cf.\ Section \ref{A2-section: additional simulations} of the supplementary material). 

\section{Discussion}\label{section: discussion}



We propose two new hybrid methods, ARGES-CIG and ARGES-skeleton, consisting of restricted versions of GES, where the restriction on the search space is not simply given by an estimated CIG or an estimated CPDAG-skeleton, but also depends adaptively on the current state of the algorithm. We include the adaptive part in our algorithms to ensure that the output is a consistent estimate of the underlying CPDAG. The fact that the adaptive modification is rather small, provides an explanation for the empirical success of inconsistent hybrid methods that restrict the search space to an estimated CIG or an estimated CPDAG-skeleton (e.g., \cite{TsamardinosEtAl06, SchulteEtAL10}).




We prove consistency of GES and ARGES in sparse high-dimensional settings with multivariate Gaussian distributions, linear SEMs with sub-Gaussian errors, or nonparanormal distributions (see Sections \ref{section: high-dimensional consistency}, \ref{section: high dimensional consistency for LSEM} and \ref{section: high dimensional consistency in the NPN setting}). To the best of our knowledge, these are the first results on high-dimensional consistency of score-based and hybrid methods. Our simulation results indicate that GES and ARGES generally outperform the PC algorithm (see Section \ref{section: high-dimensional simulations}), which has so far been the most popular structure learning algorithm in high-dimensional settings. Moreover, an advantage of (AR)GES compared to PC is that its output is always a CPDAG, which is important for some applications such as (joint-)IDA \citep{MaathuisKalischBuehlmann09, NandyMaathuisRichardson15} and the generalized adjustment criterion \citep{PerkovicEtAl15, PerkovicEtAl17}. Note that the sample version of the PC algorithm provides no such guarantee and often produces a partially directed graph that is not a CPDAG.

A disadvantage of ARGES compared to GES is that it requires an additional tuning parameter to estimate the CIG or the CPDAG-skeleton. Our simulation results suggest that ARGES-CIG can achieve a very similar performance as GES in a much shorter time as long as we choose a sufficiently dense estimated CIG (by adjusting the corresponding tuning parameter) while respecting the computation limit. In some settings, however, the restricted search space is not only beneficial from a computational point of view, but also in terms of estimation performance (see Section \ref{subsection: dense graphs}).





Tuning the penalty parameter $\lambda$ of a scoring criterion of (AR)GES is a well-known practical problem. We recommend to apply the stability selection approach of \cite{MeinshausenBuehlmann10} or 
  to use the extended BIC criterion \citep{ChenChen08, FoygelDrton10}, which has been shown to work better in sparse high-dimensional settings than the BIC criterion. 

There have been some recent theoretical and practical attempts to speed up GES and we note that they can be applied to ARGES as well. \cite{ChickeringMeek15} proposed a modification of the backward phase of GES that has polynomial complexity. Further, the authors showed that the final output of this modified version of GES, called selective GES (SGES), is consistent in the classical setting if the output of the forward phase SGES is an independence map of the CPDAG in the limit of large samples. The forward phase of ARGES can be combined with the backward phase of SGES and consistency of such an algorithm follows from the fact that the output of the forward phase of ARGES is an independence map of the CPDAG in the limit of large samples (see the proof of Theorem \ref{theorem: consistency of ARGES}). \cite{Ramsey15} showed that with an efficient implementation and parallel computing, GES can be scaled up to thousands of variables. Similar efficient implementations and parallel computations are possible for hybrid algorithms like ARGES, and this would push the computation limit even further.

 We establish a novel connection between score-based and constraint-based methods (see Section \ref{subsection: l0 penalized likelihood score}). In particular, Lemma \ref{lemma: score difference} shows that the score-based GES algorithm and the constraint-based PC algorithm are in fact closely related approaches in the multivariate Gaussian setting. The fundamental principle of the PC algorithm is to start with a complete graph and to
delete edges sequentially by testing conditional independencies. In the multivariate Gaussian setting, conditional independence tests are equivalent to tests for zero partial correlations. Lemma \ref{lemma: score difference} shows that GES checks if sample partial correlations are large enough in order to add edges in the forward phase or small enough to delete edges in the backward phase. This insight opens the door to study new score-based and hybrid algorithms that are applicable to broader classes of models. For example, in Section \ref{section: high dimensional consistency in the NPN setting}, we defined a new scoring criterion based on rank correlations (see Definition \ref{definition: NPN score}) for nonparanormal distributions. Analogously, one can define scoring criteria based on more general conditional independence tests, leading to score-based (or hybrid) competitors of the PC algorithms based on such conditional independence tests
(e.g., \cite{ZhangEtAl11, DoranEtAl14}).

Although both GES and PC use partial correlation-based conditional independence tests in the multivariate Gaussian setting, we found that GES outperforms PC in terms of estimation quality. A possible explanation for the better performance of GES is that GES considers the skeleton and the orientations of the edges simultaneously, whereas the PC algorithm first determines the skeleton and then orients the edges by determining v-structures and subsequently applying a set of rules \citep{Meek95}.



Recall that our high-dimensional consistency proofs require an assumption on the growth of oracle versions of (AR)GES (see assumption (A5) of Section \ref{subsection: high-dimensional consistency of ARGES}). Such an assumption is not required for high-dimensional consistency of the PC algorithm \citep{KalischBuehlmann07a}. We discussed this assumption in Section \ref{subsection: assumption A5} and provided some strong structural conditions under which this assumption holds. In order to derive these sufficient conditions, we show a connection between GES and the Chow-Liu algorithm, which may be of independent interest.

We emphasize that our consistency result in the classical setting (where the number of variables remains fixed and the sample size goes to infinity) does not require any distributional assumption (it even holds for discrete distributions), except that
the joint distribution of the variables is DAG-perfect, i.e., there exists a DAG $\mathcal{G}$ such that all conditional independence relationships encoded by $\mathcal{G}$ hold in the joint distribution and vice versa. \cite{ChickeringMeek02} showed that without this assumption (but with a weaker assumption on the joint distribution), GES is consistent for learning a \emph{minimal independence map} of the joint distribution. A DAG is
a minimal independence map of a distribution if it is an independence map of the distribution and no proper subgraph is an independence map. An interesting direction for future work is to investigate (in)consistency of ARGES for learning a minimal independence map under similar assumptions as in \cite{ChickeringMeek02}.

\bibliographystyle{apalike}
\bibliography{bib/Mybibliography}

\end{document}